\documentclass[11pt, a4paper]{article}

\usepackage{pdfsync}

\usepackage{amssymb,amsmath, wasysym, graphicx, epsfig, setspace, wrapfig, comment, pgfplots, theorem}
\usepackage{hyperref}
\usepackage{latexsym,dsfont,amsfonts}
\usepackage{epsfig,psfrag,color}
\usepackage{enumerate}
\usepackage[shortlabels]{enumitem}
\usepackage{algpseudocode,algorithm} 
\algdef{SE}[DOWHILE]{Do}{doWhile}{\algorithmicdo}[1]{\algorithmicwhile\ #1}

\usepackage{mathtools}

\newcommand{\abs}{\mathrm{abs}}
\newcommand{\setu}{{\mathrm{\mathfrak{u}}}}

\newcommand{\bs}[1]{\boldsymbol{#1}}

\newcommand{\supp}{\mathrm{supp}}
\newcommand{\diam}{\mathrm{diam}}
\newcommand{\tri}{\bigtriangleup}
\newcommand{\ind}{\mathrm{ind}}

\newcommand{\eps}{\epsilon}

\newcommand{\N}[0]{\mathbb{N}}

\newcommand{\R}[0]{\mathbb{R}}

\newcommand{\bszero}{{\boldsymbol{0}}}

\newcommand{\bsq}{{\boldsymbol{q}}}

\newcommand{\bsw}{{\boldsymbol{w}}}
\newcommand{\bsx}{{\boldsymbol{x}}}
\newcommand{\bsy}{{\boldsymbol{y}}}
\newcommand{\bsz}{{\boldsymbol{z}}}
\newcommand{\bsA}{{\boldsymbol{A}}}

\newcommand{\bsL}{{\boldsymbol{L}}}
\newcommand{\bsU}{{\boldsymbol{U}}}

\newcommand{\bsbeta}{{\boldsymbol{\beta}}}

\newcommand{\bsgamma}{{\boldsymbol{\gamma}}}
\newcommand{\bseta}{{\boldsymbol{\eta}}}

\newcommand{\bsnu}{{\boldsymbol{\nu}}}

\newcommand{\bsrho}{\boldsymbol{\rho}}

\newcommand{\bspsi}{\boldsymbol{\psi}}

\newcommand{\bsDelta}{{\boldsymbol{\Delta}}}

\newcommand{\bsPhi}{\boldsymbol{\Phi}}

\newcommand{\calA}{\mathcal{A}}

\newcommand{\calF}{\mathcal{F}}
\newcommand{\calG}{\mathcal{G}}
\newcommand{\calH}{\mathcal{H}}

\newcommand{\calT}{\mathcal{T}}

\newcommand{\rd}{\mathrm{d}}
\newcommand{\rmN}{\mathrm{N}}

\newcommand{\bbP}{\mathbb{P}} 
\newcommand{\bbE}{\mathbb{E}}
\newcommand{\bbR}{\mathbb{R}}

\newcommand{\ubar}[0]{\overline{u}}
\newcommand{\ellbar}[0]{\overline{\ell}}
\newcommand{\phibar}[0]{\overline{\varphi}}

\newcommand{\uhat}[0]{\widehat{u}}
\newcommand{\tauhat}[0]{\widehat{\bsq}}

\newcommand{\amin}[0]{a_{\mathrm{min}}}
\newcommand{\amax}[0]{a_{\mathrm{max}}}

\newcommand{\ellinf}[0]{\ell_{0, \mathrm{inf}}}
\newcommand{\gcdf}[0]{g_{\mathrm{cdf}}}
\newcommand{\gpdf}[0]{g_{\mathrm{pdf}}}

\definecolor{darkred}{RGB}{139,0,0}
\definecolor{darkgreen}{RGB}{0,100,0}
\definecolor{darkmagenta}{RGB}{170,0,120}
\definecolor{darkpurple}{RGB}{110,0,180}
\definecolor{darkblue}{RGB}{40,0,200}
\definecolor{darkbrown}{rgb}{0.75,0.40,0.15}

\newtheorem{theorem}{Theorem}
\newtheorem{lemma}[theorem]{Lemma}

\newtheorem{corollary}[theorem]{Corollary}

\theoremstyle{plain} \theorembodyfont{\rmfamily}
\newtheorem{assumption}{Assumption}
\newtheorem{remark}[theorem]{Remark}

\newenvironment{proof}{\begin{trivlist}
    \item[\hskip\labelsep{\bf Proof.}]}{$\hfill\Box$\end{trivlist}}

\numberwithin{theorem}{section}

\newcommand{\genvec}{\boldsymbol{\bsz}_\mathrm{gen}}

\numberwithin{equation}{section}

\title{Finite element analysis of density estimation using preintegration for elliptic PDE with random input}
\date{\today}
\author{Alexander D. Gilbert%
	\footnote{School of Mathematics and Statistics, 
			UNSW Sydney, Sydney NSW 2052, Australia.\\
                         \texttt{alexander.gilbert@unsw.edu.au}}
	     }

\begin{document}
\maketitle

\begin{abstract}
This paper analyses the finite element component of the error when
using preintegration to approximate the cdf and pdf for uncertainty quantification (UQ) problems involving elliptic PDEs with random inputs.
It is a follow up to Gilbert, Kuo, Srikumar, \emph{SIAM J.~Numer.~Anal.} \textbf{63} (2025), pp. 1025--1054, which introduced a method of density estimation for a class of UQ problems, based on computing the integral formulations of the cdf and pdf by performing an initial smoothing preintegration step and then applying a quasi-Monte Carlo
quadrature rule to approximate the remaining high-dimensional integral.
That paper focussed on the quadrature aspect of the method, whereas this paper
studies the spatial discretisation of the PDE using finite element methods.
First, it is shown that the finite element approximation satisfies the
required assumptions for the preintegration theory, including the important 
monotonicity condition. Then the finite element error is analysed and finally,
the combined finite element and quasi-Monte Carlo error is bounded.
It is shown that under similar assumptions, the cdf and pdf can be approximated
with the same rate of convergence as the much simpler problem of computing expected values.
\end{abstract}

\section{Introduction}\label{sec:Intro}
Computing the cumulative density function (cdf) and probability density function (pdf)
are important tasks in uncertainty quantification (UQ), since
the cdf and pdf provide a more complete picture of the uncertainty than 
the expected value alone.
Recently, \cite{GKSr25} introduced a method
for approximating the cdf and pdf for elliptic PDEs with random input
by using an initial smoothing by preintegration step 
followed by a quasi-Monte Carlo (QMC) rule to approximate the remaining
high-dimensional integral. This paper is a follow up to \cite{GKSr25}
and analyses the component of the approximation corresponding to the 
spatial discretisation of the PDE using finite element methods (FEM).

Let $D \subset \R^d$, for $d = 1, 2$ or $3$, be a bounded convex domain and
let $(\Omega, \calF, \pi)$ be a probability space.
Consider the random elliptic PDE
	\begin{align}
		\label{eq:pde}
		-\nabla \cdot (a(\bsx, \omega)\nabla u(\bsx, \omega)) \,&=\, \ell(\bsx, \omega),
		&\bsx \in D,
		\\\nonumber
		u(\bsx, \omega)\,&=\, 0,
		&\bsx \in \partial D,
	\end{align}
where the derivatives are with respect to the \emph{physical} variable $\bsx
\in D$, and where $\omega \in \Omega$ is a \emph{random} variable
modelling the uncertainty in the diffusion coefficient $a$ and the
source term $\ell$.
The goal of this paper is to approximate the cdf and pdf of a real-valued 
random variable representing a quantity of interest (QoI)
that is given by a linear functional $\calG$ of the solution to the PDE \eqref{eq:pde},
\begin{equation}
\label{eq:qoi}
X = \calG(u(\cdot, \omega)).
\end{equation}

The starting point for the approximation is to formulate the cdf and pdf at a point $t$
as high-dimensional integrals involving the QoI, 
\begin{align*}
F(t) &= \bbP[X \leq t] = \bbE[\ind(t - X)]
= \int_\Omega \ind(t - \calG(u(\cdot, \omega))) \, \rd \bbP(\omega),
\quad\text{and}
\\
f(t) &= \bbE[ \delta(t - X)]
= \int_\Omega \delta(t - \calG(u(\cdot, \omega))) \, \rd \bbP(\omega),
\end{align*}
where $\ind$ is the indicator function and $\delta$ is the Dirac
$\delta$ function.

The method for estimating the cdf and pdf proceeds by 
first approximating the QoI by discretising the spatial domain $D$
using finite element methods. To compute the integrals above
an initial preintegration step is performed 
to remove the discontinuities induced by the indicator and $\delta$ functions,
before finally, a quadrature rule is applied to compute the 
remaining (now smooth) high-dimensional integrals.

The focus of this paper is on the first FE component of the approximation
and since the formulations of cdf and pdf are nonsmooth functions of the PDE solution,
the preintegration step is also essential for the FEM analysis.
This paper proves that the preintegration theory
also holds for the FEM approximation of the QoI \eqref{eq:qoi} and then
bounds the FE error for approximating the cdf and pdf.
Thus, when coupled with results from \cite{GKSr25}, this paper provides a 
complete error analysis of density estimation for random elliptic PDE
\eqref{eq:pde} using preintegration along with FEM and QMC.
The main result of the paper is that despite the
added difficulty caused by the discontinuities in \eqref{eq:cdf} and \eqref{eq:pdf},
by using preintegration the cdf and pdf can be approximated with the \emph{same
rates of convergence} as for computing the expected value.

To model the case where the 
uncertainty in the diffusion coefficient and the source term are driven by two independent Gaussian random fields with (truncated) Karhunen--Lo\'eve expansions, 
we assume that the coefficient $a$ is lognormal
and $\ell$ admits a finite series expansion, each with independent parametrisations. 
Let $s \in \N$.
For a vector of independent and identically distributed (i.i.d.) standard normal 
random variables $\bsw = (w_0, w_1, \ldots, w_s) \in \R^{s + 1}$, $w_i \sim \rmN(0, 1)$,
the source term is given by
\begin{equation}
\label{eq:source}
\ell(\bsx, \omega) \equiv \ell(\bsx, \bsw) \coloneqq \ellbar(\bsx) + \sum_{i = 0}^s w_i \, \ell_i(\bsx).
\end{equation}
where $\ellbar, \ell_i \in L^2(D)$ are known and deterministic.
For a vector of independent and identically distributed (i.i.d.) standard normal 
random variables $\bsz = (z_1, z_2, \ldots, z_s) \in \R^s$, $z_i \sim \rmN(0, 1)$,
the diffusion coefficient is given by
\begin{equation}
\label{eq:coeff}
a(\bsx, \omega) \equiv a(\bsx, \bsz) \coloneqq \exp\Bigg( \sum_{j = 1}^s z_j \, a_j(\bsx) \bigg)
\end{equation}
where $a_j \in L^\infty(D)$ are known and deterministic.
To simplify the notation later, we will combine the two random vectors into
a single parameter $\bsy = (\bsw, \bsz) \in \R^{2s + 1}$.
In this way, the cdf  and pdf
are given by high-dimensional integrals with respect to the standard normal density,
$\rho(y) = 1/\sqrt{2\pi} e^{-y^2/2}$,
\begin{align}
\label{eq:cdf}
F(t) &= \int_{\R^{2s + 1}} \ind(t - \calG(u(\cdot, \bsy))) \Bigg(\prod_{i = 0}^{2s} \rho(y_i)\Bigg)
\, \rd \bsy,
\quad \text{and}
\\
\label{eq:pdf}
f(t) &= \int_{\R^{2s + 1}} \delta(t - \calG(u(\cdot, \bsy))) \Bigg(\prod_{i = 0}^{2s} \rho(y_i)\Bigg)
\, \rd \bsy.
\end{align}

Smoothing by preintegration is a technique for tackling integrals involving 
discontinuities or singularities, such as \eqref{eq:cdf} or \eqref{eq:pdf}. 
The key idea is to integrate out a single well-chosen variable, e.g., $y_0$,
in order to remove the discontinuity.
Under certain conditions, the result is a smooth function in one dimension less 
for which the integral can now
be approximated efficiently by a quadrate rule.
It is a special case of the more general conditional Monte Carlo or conditional sampling method, where one conditions on general partial information (instead of a single variable) and which has been used extensively as a method for smoothing and variance reduction
in statistics and computational finance 
\cite{ACN13a,ACN13b,BayB-HTemp23,Glasserman,GlaSta01,Hol11,LecLem00,LiuOw23,WWH17}.
Later, a sequence of papers \cite{GKS10,GKS13,GKS17note,GKLS18}
formalised the smoothing effect of preintegration.
Those papers assumed a monotone condition, i.e., that the integrand should be monotone
with respect to the preintegration variable, 
which was later shown in \cite{GKS22b} to be necessary for
the smoothing to be effective.
Recently, \cite{LecPucBAb22} introduced a practical 
method of density estimation using conditional sampling.
For the special case of preintegration, density estimation was analysed in
\cite{GKS23} and later applied to elliptic PDEs of the form \eqref{eq:pde}
in \cite{GKSr25}, which is the setting followed here.

The structure of the remainder of this paper is as follows. Section~\ref{sec:background}
summarises the key background material on PDEs with random input, finite element methods and smoothing by preintegration. Then, Section~\ref{sec:method} introduces 
the method for approximating the cdf and pdf using preintegration along with FEM in space and numerical integration.
After verifying that the FEM approximation satisfies the necessary assumptions for the preintegration theory in Section~\ref{sec:preint-theory},
Section~\ref{sec:error} provides a full analysis of the error, including the FE and combined FE-QMC error.

\section{Background on random PDEs and smoothing by preintegration}
\label{sec:background}

\subsection{PDEs with random input}

This section introduces the necessary background material and assumptions on the parametric PDE \eqref{eq:pde}.
To ensure that the weak form of \eqref{eq:pde} is well posed and that 
the solution is sufficiently regular for the analysis later on,
we make the following assumptions on the coefficient $a$, the source term $\ell$ 
and the physical domain $D$. In what follows let $V \coloneqq H^1_0(D)$, the first-order Sobolev space 
with square-integrable derivatives and vanishing boundary trace, and let $V^* $ denote the dual space. 
The duality pairing on $V^* \times V$ is denoted by $\langle \cdot, \cdot \rangle$ and is continuously extended
from the usual $L^2(D)$ inner product, which, with a slight abuse of notation, 
is also denoted by $\langle \cdot, \cdot \rangle$.
Also, let $W^{1, \infty}(D)$ denote the
first-order Sobolev space with essentially bounded weak derivatives, 
with norm given by $\|v\|_{W^{1, \infty}(D)} \coloneqq \max \{ \|v\|_{L^\infty(D)}, \|\nabla \|_{L^\infty(D)}\}$.

\begin{assumption}
\label{asm:pde} Let $d= 1, 2, $ or $3$ and assume that:
\begin{enumerate}[(a)]
\item $D \subset \R^d$ is a bounded, strictly convex polygonal or polyhedral domain;
\item for $\bsz \in \R^s$, $a(\cdot, \bsz) \in W^{1, \infty}(D)$ is of the form \eqref{eq:coeff} with $a_j \in W^{1, \infty}(D)$
and $z_j \sim \mathrm{N}(0, 1)$ for $j = 1, 2, \ldots, s$;
\item for $\bsw \in \R^{s + 1}$, $\ell(\cdot, \bsw) \in V^*$ is of the form \eqref{eq:source}
with $\ellbar, \ell_i \in V^*$ and $w_i \sim \mathrm{N}(0, 1)$ for $i = 0, 1, 2, \ldots, s$; and
\item $\ell_0 \in L^\infty(D)$ with $\ell_{0, \inf} \coloneqq \inf_{\bsx \in D} \ell_0(\bsx) > 0$.
\end{enumerate}
\end{assumption}

An important consequence of Assumption~\ref{asm:pde}(b), is that for
$\bsz \in \R^s$ the coefficient $a(\cdot, \bsz)$ is bounded from above and below,
i.e., there exist functions $\amin, \amax $ such that
\begin{equation}
\label{eq:a_bound}
0 \leq \amin(\bsz) \leq a(\bsx, \bsz) \leq \amax(\bsz) < \infty 
\quad \text{for all } \bsx \in D.
\end{equation}
In addition, the bounds satisfy $1/\amin, \amax \in L^p_{\bsrho}(\R^s)$ for all $1 \leq p < \infty$, where
$L^p_{\bsrho}(\R^s)$ is the space of $p$-integrable functions with respect to the product normal density $\bsrho$.

Multiplying the PDE \eqref{eq:pde} by $v \in V$ then integrating over $D$, using integration by parts, 
gives the variational equation
\[
\int_D a(\bsx, \bsz) \nabla u(\bsx, \bsy) \cdot \nabla v(\bsx) \, \rd \bsx = 
\int_D \ell(\bsx, \bsw) \, v(\bsx) \, \rd \bsx,
\]
where recall that we define combined parameter $\bsy \in \R^{2s + 1}$ by
\[
\bsy = (y_0, y_1, \ldots, y_{2s}) = (w_0, w_1, \ldots, w_s, z_1, \ldots, z_s) = (\bsw, \bsz).
\]
For $\bsz \in \R^s$, the parametric bilinear form $\calA(\bsz; \cdot, \cdot) : V \times V \to \R$ is defined by
\begin{equation}
\label{eq:bilin}
\calA(\bsz; v, v') \, \coloneqq\, \int_D a(\bsx, \bsz) \nabla v(\bsx) \cdot \nabla v'(\bsx) \, \rd \bsx,
\end{equation}
which by \eqref{eq:a_bound} is bounded and coercive for each $\bsz \in \R^s$, i.e., for 
any $v, v' \in V$
\begin{equation}
\label{eq:bound+coerc}
|\calA(\bsz; v, v')| \leq \amax(\bsz) \|v\|_V \|v'\|_V
\quad \text{and} \quad
\calA(\bsz; v, v) \geq \amin(\bsz) \|v\|_V^2.
\end{equation}

For each $\bsy \in \R^{2s + 1}$, the weak form of \eqref{eq:pde} is:
find $u(\cdot, \bsy) \in V$ such that
\begin{equation}
\label{eq:varpde}
\calA(\bsz; u(\cdot, \bsy), v) = \langle \ell(\cdot, \bsw), v\rangle
\quad \text{for all } v \in V.
\end{equation}
Since $\calA(\bsz; \cdot, \cdot)$ is bounded and coercive \eqref{eq:bound+coerc},
for each $\bsy \in \R^{2s + 1}$, the weak equation \eqref{eq:varpde} above admits a unique
solution $u(\cdot, \bsy) \in V$, which satisfies the \emph{a priori} bound
\begin{equation}
\label{eq:LaxMil}
\|u(\cdot, \bsy)\|_V \leq \frac{\|\ell(\cdot, \bsw)\|_{V^*}}{\amin}
\leq \frac{1}{\amin(\bsz)} \Bigg( \| \ellbar\|_{V^*} + \sum_{i = 0}^s |w_i| \|\ell_i\|_{V^*}\Bigg),
\end{equation}
where the second inequality follows by the linearity of the source term $\ell(\cdot, \bsw)$
\eqref{eq:source}. Similarly, for $\calG \in V^*$ the QoI of the solution satisfies
\begin{equation}
\label{eq:LaxMil-G}
|\calG(u(\cdot, \bsy))| \leq \|\calG\|_{V^*} \| u(\cdot, \bsy)\|_{V} 
\leq \frac{\|\calG\|_{V^*} \|\ell(\cdot, \bsw)\|_{V^*}}{\amin}.
\end{equation}

Since in Assumption~\ref{asm:pde} it is assumed that $D$ is convex and
$a(\cdot, \bsy) \in W^{1, \infty}(D)$, the smoothness of the solution (up to order 2) will
be governed by the regularity of the source term~$\ell(\cdot,\bsw)$.
The following parameter-dependent regularity result is classical and can
be found in, e.g, \cite{GilbTrud}. The explicit dependence on the
stochastic parameter $\bsy$ is given by a simplification of
\cite[Theorem~2.1]{TSGU13}. To allow for fractional smoothness we define 
the following spaces. For $\tau \in (0, 1)$, $C^{0, \tau}(D)$ denotes the space of H\"older
continuous functions with exponent $\tau$, and for $k \in \N_0$, $H^{k + \tau}(D)$ denotes the
usual fractional order Sobolev space (for $p = 2$) with the dual space 
denoted by $H^{-(k + \tau)}(D) = (H^{k + \tau}(D))^*$.

\begin{theorem}
\label{thm:u-spatial} Suppose that Assumption~\ref{asm:pde} holds and for
$\bsy = (\bsw, \bsz) \in \R^{2s + 1}$ let $\ell(\cdot, \bsw) \in H^{-1 +
\tau'}(D)$ for some $\tau' \in (0, 1]$. Then, the solution of
\eqref{eq:varpde} satisfies $u(\cdot, \bsy) \in V \cap H^{1 +
\tau'}(D)$ with
\begin{equation}
\label{eq:u-spatial}
\|u(\cdot, \bsy)\|_{H^{1 + \tau'}(D)} \leq C_{\tau'}
\frac{\amax(\bsz)\, \|a(\cdot,\bsz)\|_{C^{0, \tau'}(\overline{D})}^2}{\amin^4(\bsz)}\,
\|\ell(\cdot, \bsw)\|_{H^{-1 + \tau'}(D)},
\end{equation}
where $C_{\tau'} < \infty$ is independent of $\bsy$.
\end{theorem}

By the linearity of the source term \eqref{eq:source}, we can write the 
solution $u(\cdot, \bsy)$ as a sum of $s + 2$ solutions,
\begin{equation}
\label{eq:u-linear}
u(\bsx, \bsy) = \ubar(\bsx, \bsz) + \sum_{i = 0}^s w_i u_i(\bsx, \bsz)\,
\end{equation}
where $\ubar(\cdot, \bsz) \in V$ and $u_i(\cdot, \bsz)\in V$ correspond to solutions of \eqref{eq:varpde}
with deterministic  source terms $\ellbar \in V^*$ and $\ell_i \in V^*$, respectively.
Note that $\ubar(\cdot, \bsz)$ and $u_i(\cdot, \bsz)$ are independent of $\bsw$ but still depend on the parameter 
$\bsz$ via the coefficient $a(\cdot, \bsz)$. 
All of the \emph{a priori} and regularity results presented in this section also hold for
the solutions $\ubar(\cdot, \bsz) $ and $u_i(\cdot, \bsz)$ by simply replacing $\ell(\cdot, \bsw)$
by $\ellbar$ or $\ell_i$ as appropriate.

\subsection{Finite element methods}
\label{sec:fem}

To compute the solution of the PDE \eqref{eq:pde}, we discretise the physical domain
$D$ and then solve the PDE numerically using the finite element method (FEM).
Let $\{V_h\}_{h > 0}$ be a family of conforming finite element (FE)
spaces, where each space has dimension $M_h < \infty$ and corresponds to a
shape regular triangulation $\calT_h$ of the spatial domain~$D$. The
parameter $h$ is called the \emph{meshwidth}, given by $h = \max\{
\diam(\tri) : \tri \in \calT_h\}$. In this paper we consider the
continuous, piecewise linear FEM, i.e., each $V_h$ is the space of
continuous functions on $D$ that are linear on each element $\tri \in
\calT_h$. The preintegration method outlined in this paper works
for more general spaces $V_h$, e.g., higher-order polynomials, 
however, to obtain faster convergence rates
higher smoothness assumptions are required on the domain $D$ as well as on the coefficient $a(\cdot, \bsy)$ and 
the source term $\ell(\cdot, \bsy)$.

For $\bsy = (\bsw, \bsz)  \in \R^{2s + 1}$, the FE approximation of
\eqref{eq:varpde} is: find $u_h(\cdot, \bsy) \in V_h$ such that
\begin{equation}
\label{eq:fem}
\calA(\bsz; u_h(\cdot, \bsy), v_h) = \langle \ell(\cdot, \bsw), v_h \rangle
\quad
\text{for all } v_h \in V_h.
\end{equation}

This leads to a linear system defined as follows.
For $h > 0$, denote the FE basis functions by $\phi_{k} \coloneqq \phi_{h, k} \in V_h$
for $k = 1, 2, \ldots, M_h$, which correspond to the hat functions
on the FE nodes in the interior of $D$.
The hat functions on the boundary nodes are similarly denoted by 
$\phi_{k} \coloneqq \phi_{h, k}$
for $k = M_h + 1, \ldots, \overline{M_h}$. 
The discrete problem \eqref{eq:fem} is equivalent to the linear system
\begin{equation}
\label{eq:fem-matrix}
\bsA_h(\bsz) \bsU_h(\bsy) = \bsL_h(\bsw)
\end{equation}
where, recalling that $\bsy = (\bsw, \bsz)$, we define $\bsU_h(\bsy) = [U_{h, k}(\bsy)]_{k = 1}^{M_h}$ to be the coefficients
for the interior nodes,
\begin{align}
\label{eq:stiff}
\bsA_h(\bsz) \,&\coloneqq\, [A_{k, m}(\bsz)]_{k, m = 1}^{M_h}
\quad &&\text{with} \quad A_{k, m}(\bsz) \,\coloneqq\, \calA(\bsz; \phi_k, \phi_m),
\quad \text{and}
\\
\nonumber\bsL_h(\bsy) \,&\coloneqq\, [L_k(\bsy)]_{k = 1}^{M_h}
\quad &&\text{with}\quad
L_k(\bsw) \,\coloneqq\, \langle \ell(\cdot, \bsw), \phi_k \rangle .
\end{align}
The solution $u_h(\bsx, \bsy)$ of \eqref{eq:fem} is given by
\begin{equation}
\label{eq:u_h}
u_h(\bsx, \bsy) = \sum_{k = 1}^{\overline{M_h}} U_{h, k}(\bsy) \phi_k(\bsx), 
\end{equation}
where since \eqref{eq:pde} has zero Dirichlet boundary conditions,
the coefficients corresponding to the boundary nodes are $U_{h, k}(\bsy) = 0$
for $k = M_h + 1, \ldots, \overline{M_h}$.
Similar to the continuous problem, since $\calA(\bsz;\cdot,\cdot)$ is
bounded and coercive \eqref{eq:bound+coerc}, for each $\bsy = (\bsw, \bsz)
\in \R^{2s + 1}$, the FE problem \eqref{eq:fem} admits a unique solution
$u_h(\cdot, \bsy) \in V_h$ satisfying
\begin{equation}
\label{eq:LaxMil-fe}
\|u_h(\cdot, \bsy)\|_V \leq \frac{1}{\amin(\bsz)} \|\ell(\cdot,\bsw)\|_{V^*}.
\end{equation}

Classic convergence results from FE analysis give error estimates for the solutions
of \eqref{eq:fem} for a fixed stochastic parameter $\bsy$.
Then, to account for the stochastic dependence, convergence results in $L^p$
spaces have been given in \cite{TSGU13}.
To facilitate our error analysis of the FE approximations of the cdf \eqref{eq:cdf} and pdf \eqref{eq:pdf}
later, we require the following FE convergence results with explicit dependence on $\bsy$.
The following theorem follows from intermediate results in \cite{TSGU13}.

\begin{theorem}
\label{thm:fe_y} Suppose
that Assumption~\ref{asm:pde} holds and for all $\bsw \in \R^{s + 1}$
assume $\ell(\cdot, \bsw) \in H^{-1 + \tau'}(D)$ for some $\tau' \in (0, 1]$. 
Then the solution $u_h(\cdot, \bsy) \in V_h$ to the FE problem \eqref{eq:fem} satisfies
\begin{equation}
\label{eq:fe-err-u}
\bbE \big[ \|u - u_h\|_V\big]
= \int_{\bbR^{2s+1}}  \|u(\cdot,\bsy) - u_h(\cdot,\bsy)\|_V\, \Bigg(\prod_{i=0}^{2s} \rho(y_i)\Bigg) \rd\bsy
\leq C_u\, h^{\tau'},
\end{equation}
where $C_u < \infty$ is independent of $\bsy$ and $h$.
Additionally, for $\calG \in H^{-1 + \tau}(D)$ for some $\tau \in [0, 1]$
\begin{equation}
\label{eq:fe-err-qoi-point}
\big|\calG(u(\cdot, \bsy)) - \calG(u_h(\cdot, \bsy))\big|
\leq
C_\calG \frac{\amax^4(\bsz) \|a(\cdot, \bsz)\|_{C^{0, 1}(\overline{D})}^4}{\amin^9(\bsz)}
\|\ell(\cdot,\bsw)\|_{H^{-1 + \tau'}(D)} 
h^{\tau' + \tau},
\end{equation}
where $C_\calG < \infty$ is independent of $\bsy$ and $h$.
\end{theorem}
\begin{proof}
The first result \eqref{eq:fe-err-u} is given by \cite[Theorem~2.2]{TSGU13}
with $p = 1$.

The second result follows by the classic Aubin--Nitsche duality argument,
given for this setting in \cite[Section~3.1]{TSGU13}.
For brevity, we do not give the full argument here, but only quote
relevant results from \cite{TSGU13}.
Theorem~\ref{thm:u-spatial} implies that $u(\cdot, \bsy) \in H^{1 + \tau'}(D)$,
and so using again \cite[Theorem~2.2]{TSGU13} we have
\begin{align}
\label{eq:u-fe-err-point}
\|u(\cdot, \bsy) - u_h(\cdot, \bsy)\|_V
\,&\lesssim\, \bigg(\frac{\amax(\bsz)}{\amin(\bsz)}\bigg)^{1/2}
\|u(\cdot, \bsy)\|_{H^{1 + \tau'}(D)} \, h^{\tau'}
\nonumber\\
&\lesssim\, \bigg(\frac{\amax^3(\bsz)}{\amin^9(\bsz)}\bigg)^{1/2}
\|a(\cdot,\bsz)\|_{C^{0, 1}(\overline{D})}^2\, \|\ell(\cdot, \bsw)\|_{H^{-1 + \tau'}(D)} \, h^{\tau'},
\end{align}
where each $\lesssim$ indicates that we have suppressed a constant
factor which is independent of both $\bsy$ and $h$. In the second step we
used \eqref{eq:u-spatial} together with $\|a(\cdot, \bsz)\|_{C^{0,
\tau'}(\overline{D})} \leq \diam(D)^{1 - \tau'}\|a(\cdot, \bsz) \|_{C^{0,
1}(\overline{D})}$ for $\tau' \leq 1$.

Similarly, since $\calA(\bsy; \cdot, \cdot)$ is symmetric,
the solution $u^\calG(\cdot, \bsy) \in V$ to the dual
problem,
\[
\calA(\bsy; v, u^\calG(\cdot, \bsy)) = \calG(v) \quad \text{for all }
v \in V,
\]
and its FE approximation $u^\calG_h(\cdot, \bsy) \in V_h$ also satisfy
\eqref{eq:u-fe-err-point}, but with $\ell(\cdot, \bsw)$ replaced by $\calG \in H^{-1 + \tau}(D)$.

Lemma 3.2 from \cite{TSGU13} gives the error bound
\[
|\calG(u(\cdot, \bsy)) - \calG(u_h(\cdot, \bsy))|
\leq \amax(\bsz) \|u(\cdot, \bsy) - u_h(\cdot, \bsy)\|_V
\|u^\calG(\cdot, \bsy) - u^\calG_h(\cdot, \bsy)\|_V,
\]
from which the result \eqref{eq:fe-err-qoi-point} follows after
using \eqref{eq:u-fe-err-point} for both
the primal and dual FE errors, with exponents $\tau'$ and $\tau$, respectively.
\end{proof}

\begin{remark}
Note that this theorem also holds for deterministic right-hand sides by replacing $\ell(\cdot, \bsw)$
by one of $\ellbar$ or $ \ell_i$, $ i = 0, 1, \ldots, s$.
\end{remark}

\subsection{Density estimation using smoothing by preintegration} 
\label{sec:preint}
Here we briefly summarise the strategy from \cite{GKS23} for approximating
the cdf and pdf of a random variable, given by $X = \varphi(y_0, \bsy)$, 
where for convenience we now separate out the first variable $y_0$ and 
denote $\bsy = (y_1, y_2, \ldots, y_{2s})$. Let $[t_0, t_1]$ be 
a finite interval on which we are interested computing the density.
The cdf and pdf at a point $t \in [t_0, t_1]$ are obtained by computing the integrals:
\begin{align}
\label{eq:cdf0}
F(t) &= \int_{\R^{2s + 1}} \ind\big(t - \varphi(y_0, \bsy)\big) \Bigg( \prod_{ i = 0}^{2s} \rho(y_i)\Bigg) \, \rd y_0 \rd \bsy,
\\
\label{eq:pdf0}
f(t) &= \int_{\R^{2s + 1}} \delta\big(t - \varphi(y_0, \bsy)\big) \Bigg( \prod_{ i = 0}^{2s} \rho(y_i)\Bigg) \, \rd y_0 \rd \bsy.
\end{align}
To begin with, we consider a generic function $\varphi$, but later we will
consider the specific case of a QoI given by a functional of the discrete approximation.

The key idea behind preintegration is to first integrate the \emph{nonsmooth} integrands above
with respect to a single variable, taken here without loss of generality to be $y_0$, such
that the results are \emph{smooth} functions of the remaining $2s$ variables $\bsy = (y_1, y_2, \ldots, y_{2s})$.
The preintegrated functions are defined as
\begin{align}
\label{eq:gcdf0}
\gcdf(\bsy) &\coloneqq\, \int_{-\infty}^\infty \ind\big(t - \varphi(y_0, \bsy)\big) \rho(y_0) \, \rd y_0,
\quad \text{and}
\\
\label{eq:gpdf0}
\gpdf(\bsy) &= \int_{-\infty}^\infty \delta\big(t - \varphi(y_0, \bsy)\big) \rho(y_0) \, \rd y_0,
\end{align}
so that the cdf \eqref{eq:cdf0} and pdf \eqref{eq:pdf0} are 2s-dimensional integrals
of smooth functions
\begin{equation}
F(t) = \int_{\R^{2s}} \gcdf(\bsy) \Bigg(\prod_{i = 1}^{2s} \rho(y_i)\Bigg) \rd \bsy
\quad \text{and} \quad
f(t) = \int_{\R^{2s}} \gpdf(\bsy) \Bigg(\prod_{i = 1}^{2s} \rho(y_i)\Bigg) \rd \bsy.
\end{equation}
Since after preintegration the resulting functions are smooth, the integrals above 
can be approximated efficiently by a quadrature rule, such as QMC \cite{GKS23,GKSr25}
or sparse grids \cite{BG04,Hol11}.

To characterise the regularity with respect to the parameter $\bsy$, we define the
following notation and function spaces. Let $\N_0 = \{0, 1, 2, \ldots\}$. For $i = 0, 1, \ldots, 2s$ and a multiindex $\bsnu \in \N_0^{2s + 1}$, denote the first-order and order $\bsnu$ mixed derivatives by
\[
\partial^i \coloneqq \frac{\partial }{\partial y_i}
\quad \text{and} \quad 
\partial ^\bsnu \coloneqq \prod_{i = 0}^{2s} \frac{\partial^{\nu_i}}{\partial y_i^{\nu_i}},
\] 
respectively. The weak derivative of $\phi : \R^{2s + 1} \to \R$ is the distribution
satisfying
\[
\int_{\R^{2s + 1}} \partial^\bsnu \phi(\bsy) v(\bsy) \, \rd \bsy 
= (-1)^{|\bsnu|} \int_{\R^{2s + 1}} \phi(\bsy) \partial ^\bsnu v(\bsy)\, \rd \bsy
\quad \text{for all } v \in C^\infty_0(\R^{2s + 1}),
\]
where $C^\infty_0(\R^{2s + 1})$ is the space of $\infty$-differentiable functions with compact support.

For $\bsnu \in \N_0^{2s + 1}$, let $C^\bsnu(\R^{2s + 1})$ denote the space of
functions with continuous mixed (classical) derivatives up to order $\bsnu$,
and let $\calH_{2s + 1}^\bsnu$ denote the \emph{weighted Sobolev space
of dominating mixed smoothness} of order $\bsnu$. The weighted
Sobolev space $\calH_{2s + 1}^\bsnu$ depends on a collection of
weight parameters $\bsgamma \coloneqq \{\gamma_\setu > 0 : \setu \subseteq \{1:2s\}\}$,
which model the importance of different subsets of variables,
and a weight function $\psi : \R^{2s + 1} \to (0, \infty)$, which is strictly positive and controls
the growth of derivatives as $y_i \to \pm \infty$. 
Formally, $\calH_{2s + 1}^\bsnu$ is defined to be the space of locally integrable functions on $\R^{2s + 1}$ such that
the norm
\begin{equation}
\| \phi\|_{\calH_{2s + 1}^\bsnu}^2 \coloneqq 
\sum_{\bseta \leq \bsnu} \frac{1}{\gamma_\bseta}
\int_{\R^{2s + 1}}| \partial ^ \bseta \phi(\bsy)|^2
\bspsi_\bseta(\bsy_\bseta) \bsrho_{-\bseta}(\bsy_{-\bseta}) \, \rd \bsy,
\end{equation}
is finite, where we define notation
\[
\gamma_\bseta \coloneqq \gamma_{\supp(\bseta)},\quad
\bspsi_\bseta(\bsy_\bseta) = \prod_{i = 0, \eta_i \neq 0}^{2s} \psi(y_i),
\quad \text{and} \quad
\bsrho_{-\bseta}(\bsy_{-\bseta}) \coloneqq \prod_{i = 0, \eta_i = 0}^{2s} \rho(y_i).
\]
In this paper, we take a Gaussian weight function $\psi(y) \coloneqq \exp(-\mu y^2)$ for some $\mu \in (0, 1/2)$,
although other choices are possible, e.g., exponential weight functions, see \cite{KSWW10}.
The $2s$-dimensional spaces are defined analogously by omitting $y_0$. 

The following assumptions on $\varphi$ ensure that the preintegration step is well defined and can be performed in practice.

\begin{assumption}
\label{asm:preint}
For $s \geq 1$ and $\bsnu \in \N_0^{2s}$, assume that $\varphi : \R^{2s + 1} \to \R$ satisfies:
\begin{enumerate}[(a)]
\item $\partial^0 \varphi(y_0, \bsy) > 0$ for all $(y_0, \bsy) \in \R^{2s + 1}$;

\item for each $\bsy \in \R^{2s}$, $\varphi(y_0, \bsy) \to \infty$ as $y_0 \to \infty$; and

\item $\varphi \in \calH_{2s + 1}^{(\nu_0, \bsnu)} \cap C^{(\nu_0, \bsnu)}(\R^{2s + 1})$ for $\nu_0 = |\bsnu| + 1$.
\end{enumerate}

\end{assumption}

The monotone assumption above implies that along $y_0$  the point of where the discontinuity in the integrands in
\eqref{eq:cdf0} and \eqref{eq:pdf0} occurs is either unique or does not exist. For $t \in [t_0, t_1]$, the set of $\bsy \in \R^{2s}$ such that the discontinuity occurs is denoted
$U_t \coloneqq \{ \bsy \in \R^{2s} : \varphi(y_0, \bsy) = t $ for some $y_0 \in \R\}$,
and we also define the set of all pairs $(t, \bsy)$ such that the discontinuity occurs by
$R = \{ (t, \bsy) \in (t_0, t_1) \times \R^{2s} : \varphi(y_0, \bsy) = t $ for some $y_0 \in \R \}$.
The implicit function theorem \cite[Theorem 3.1]{GKS23} (which is an extension of \cite[Theorem 2.3]{GKS10})
implies that if $R \neq \emptyset$ then there exists a unique function $\xi \in C^{(\nu_0, \bsnu)}(\overline{R})$
defining the point of discontinuity $y_ 0 = \xi(t, \bsy)$, i.e., $\xi$ satisfies
\begin{equation}
\label{eq:xi}
\varphi(\xi(t, \bsy), \bsy)) = t
\quad \text{for } (t, \bsy) \in \overline{R}.
\end{equation}

It follows that the preintegrated functions in \eqref{eq:gcdf0} and \eqref{eq:gpdf0} 
are given by
\begin{equation}
\label{eq:g_preint}
\gcdf(\bsy) = \Phi(\xi(t, \bsy))
\quad \text{and} \quad 
\gpdf(\bsy) = \frac{\rho(\xi(t, \bsy))}{\partial^0 \varphi(\xi(t, \bsy), \bsy))}
\end{equation}
if $\bsy \in U_t$ and $0$ otherwise. It can be shown that
these preintegrated functions are as smooth as the original $\varphi$ but in one dimension less,
i.e., if Assumption~\ref{asm:preint} holds for $\bsnu \in \N_0^{2s}$, then 
$\gcdf \in \calH_{2s}^\bsnu \cap C^\bsnu(\R^{2s})$ and $\gpdf \in \calH_{2s}^\bsnu \cap C^\bsnu(\R^{2s})$,
see \cite{GKLS18} and \cite{GKS23}, respectively.

As we will see in the following section, for our problem the discontinuity $\xi$ can be computed explicitly.
However, in general no closed form is available, in which case \eqref{eq:xi} can be solved numerically by, e.g., Newton's method, as in \cite{BayB-HTemp23,GKLS18}.

\section{Density estimation for the discrete problem}
\label{sec:method}
Now we outline how to use preintegration to approximate the cdf and pdf of a
random variable given by the FE approximation of the QoI \eqref{eq:qoi},
\begin{equation}
\label{eq:X_h}
X_h =  \varphi_h(y_0, \bsy) \coloneqq \calG(u_h(\cdot, (y_0, \bsy))) ,
\end{equation}
for $\calG \in V^*$.
After applying preintegration as described in Section~\ref{sec:preint},
the FE approximations of the cdf and pdf are defined by
\begin{align}
\label{eq:F_h-preint}
F(t) \,&\approx\, F_h(t)
\,\coloneqq\, \int_{\R^{2s}} g_\mathrm{cdf}^h(\bsz) \bsrho(\bsy) \, \rd \bsy, 
\\
\label{eq:f_h-preint}
f(t) \,&\approx\; f_h(t)
\,\coloneqq\,  \int_{\R^{2s}} g_\mathrm{pdf}^h(\bsz) \bsrho(\bsy) \, \rd \bsy,
\end{align}
where $\gcdf^h$ and $\gpdf^h$ are the FEM approximations of
the preintegrated functions in \eqref{eq:g_preint}.
Analogous to the continuous problem from \cite{GKSr25}, 
using the linear structure of the source term in \eqref{eq:fem} we can simplify the preintegration step.

Similar to \eqref{eq:u-linear}, 
the FEM solution can also be written as a sum of solutions
\begin{equation}
\label{eq:u_h-linear}
u_h(\cdot, (y_0, \bsy)) = \ubar_h(\bsx, \bsz) + \sum_{i = 0}^s w_i u_{ i, h}(\bsx, \bsz),
\end{equation}
where $\ubar_h(\cdot, \bsz) \in V_h$ and $u_{h, i}(\cdot, \bsz)\in V_h$ correspond to solutions of FEM problem \eqref{eq:fem}
with deterministic source terms $\ellbar \in V^*$ and $\ell_i \in V^*$, respectively.
As before, $\ubar_h(\cdot, \bsz)$ and $u_{i, h}(\cdot, \bsz)$ are independent of $\bsw$ but still depend on the parameter 
$\bsz$ via the coefficient $a(\cdot, \bsz)$.
Since $\calG$ is linear, the discrete QoI \eqref{eq:X_h} can be expanded in a similar way,
\begin{equation}
\label{eq:phi_h-linear}
\varphi_h(y_0, \bsy) = \phibar_h(\bsz) + \sum_{i = 0}^s w_i \varphi_{i, h}(\bsz),
\end{equation}
where $\phibar_h(\bsz) \coloneqq \calG(\ubar_h(\cdot, \bsz))$ and $\varphi_{i, h}(\bsz) \coloneqq \calG(u_{i, h}(\cdot, \bsz))$. Note that the FE error bound \eqref{eq:fe-err-qoi-point} also holds for
$\phibar_h(\bsz) = \calG(\ubar_h(\cdot, \bsz))$ and $\varphi_{i, h}(\bsz) = \calG(u_{i, h}(\cdot, \bsz))$
by replacing $\ell(\cdot, \bsw)$ by either $\ellbar$ or $\ell_i$.

Next, the FEM approximations of the preintegrated functions in the cdf \eqref{eq:F_h-preint} and pdf \eqref{eq:f_h-preint}
are given by
\begin{equation}
\label{eq:g_h-preint}
g_\mathrm{cdf}^h(\bsz) \,\coloneqq\, \Phi(\xi_h(t, \bsy))
\quad \text{and} \quad
 g_\mathrm{pdf}^h(\bsz) \,\coloneqq\,
 \frac{\rho(\xi_h(t, \bsy))}{\varphi_{0,h}(\bsz)},
\end{equation}
where $\xi_h$ denotes the unique point of discontinuity in the $y_0$ direction
for the discrete QoI \eqref{eq:X_h}. Note that in \eqref{eq:g_h-preint} above,
we have simplified $\gpdf^h$ compared to the general form in \eqref{eq:g_preint} as follows. 
Recalling that $y_0 = w_0$, it follows from \eqref{eq:phi_h-linear} that 
\[
\partial^0 \varphi(y_0, \bsy) = \frac{\partial}{\partial w_0} \Bigg(\phibar_h(\bsz) + \sum_{i = 0}^s w_i \varphi_{i, h}(\bsz)\Bigg)
= \varphi_{0, h}(\bsz).
\]

Using the expansion \eqref{eq:phi_h-linear}, for each $t \in [t_0, t_1]$ and $\bsy \in \R^{2s}$,
 we can find the point of discontinuity $\xi_h(t, \bsy)$ by solving $\varphi_h(\xi_h(t, \bsy), \bsy) = t$ directly,
\begin{equation}
\label{eq:xi_h}
\xi_h(t, \bsy)
\,\coloneqq\, \frac{t-\phibar_h(\bsz)-\sum_{i=1}^s w_i\,\varphi_{i, h}(\bsz)}{\varphi_{0, h}(\bsz)}.
\end{equation}
In order for $\xi_h$ to be well defined,
we require that $\varphi_{0, h}(\bsz) > 0$ for all $\bsz$, which will be verified in the following section.
Note also that in this case for each $t \in [t_0, t_1]$ there is always a value $\bsy \in \R^{2s}$
such that the discontinuity occurs, i.e., $U_t = \R^{2s}$.
Therefore, for all $t \in [t_0, t_1]$ and $\bsy \in \R^{2s}$ the preintegrated functions in $\gcdf^h(\bsy)$ and $\gpdf^h(\bsy)$ in 
\eqref{eq:g_h-preint} are well defined and can be computed explicitly. 
Note that the simplifications above also hold for the continuous QoI $X = \varphi (\bsy)$ \eqref{eq:qoi},
again see \cite{GKSr25}.

The final step in approximating the cdf and pdf is to compute numerically the remaining integrals in 
\eqref{eq:F_h-preint} and \eqref{eq:f_h-preint} using a quadrature rule. 
Let $Q_{2s, N}$ be some quadrature rule using $N$ points designed to compute $2s$-dimensional integrals 
with respect to the standard normal density, e.g., a randomly shifted lattice QMC rule, as used in
\cite{GKSr25}, given by
\begin{equation}
\label{eq:QMC}
Q_{2s, N} (g) = \frac{1}{N} \sum_{k = 0}^{N - 1} g ( \bsq_k),
\quad \text{where } \bsq_k \coloneqq \bsPhi^{-1}\bigg(\bigg\{\frac{k \genvec}{N} + \bsDelta\bigg\}\bigg),
\end{equation}
for a \emph{generating vector} $\genvec \in \{1, 2, \ldots, N - 1\}^{2s}$ with components coprime to $N$ and a uniformly distributed \emph{random shift} $\bsDelta \in [0, 1)^{2s}$.
Here $\bsPhi^{-1}$ is the standard normal inverse cdf (used to transform QMC points the unit cube to $\R^{2s}$)
and $\{\cdot\}$ denotes taking the fractional part, both applied componentwise to a vector.
See, e.g., \cite{DKS13,DP10} for further details. Other potential choices for quadrature include sparse grid methods (e.g., based on 1D Gauss--Hermite rules).

Then the full approximations of the cdf and pdf at $t \in [t_0,t_1]$
using a FE with meshwidth $h > 0$ and $N \in \N$ quadrature points are
respectively given by
\begin{align}
\label{eq:F_hN}
F(t) \,&\approx\, F_{h, N}(t)
\,\coloneqq\, Q_{2s, N} (g_\mathrm{cdf}^h)
= Q_{2s, N} \big(\Phi(\xi_h(t, \cdot))\big),
\\
\label{eq:f_hN}
f(t) \,&\approx\; f_{h, N}(t)
\,\coloneqq\, Q_{2s, N} (g_\mathrm{pdf}^h)
= Q_{2s, N}
\bigg( \frac{\rho(\xi_h(t, \cdot))}{\varphi_{0,h}(\cdot)}\bigg).
\end{align}

\section{Preintegration theory for the discrete problem}
\label{sec:preint-theory}

In this section, we verify that the discrete QoI \eqref{eq:X_h} satisfies Assumption~\ref{asm:preint},
which will ensure that the preintegration step outlined in the previous section is well defined.

First, we show that the QoI \eqref{eq:X_h} is monotone with respect to $y_0$, i.e.,
we prove that $\partial^0 \varphi_h(y_0, \bsy) = \varphi_{0, h}(\bsz) > 0$ for all $\bsz \in \R^s$.
Throughout this section we will also utilise the corresponding result for
the continuous solution.
In \cite[Lemma 5.3 and Corollary 5.4]{GKSr25}, it was shown that
for a positive linear functional $\calG \in V^*$ 
the solution $u_0(\cdot, \bsz) \in V$ to \eqref{eq:varpde} with
right-hand side $\ell_0$ satisfies the lower bound
\begin{equation}
\label{eq:phi0-lower}
\varphi_0(\bsz) \coloneqq \calG(u_0(\cdot, \bsz)) \geq \frac{\calG(u_0(\cdot, \bszero))}{K_0(\bsz)}
\quad \text{for all } \bsz \in \R^{s},
\end{equation}
where
\begin{align} \label{eq:K0}
K_0(\bsz) \,\coloneqq\, \amax(\bsz)\bigg( 1
+ \frac{\|u_0(\cdot, \bszero)\|_{W^{1, \infty}(D)}}{\ell_{0, \mathrm{inf}}}
\sum_{j = 1}^s |z_j|\, \|a_j\|_{W^{1, \infty}(D)}
\bigg).
\end{align}

\subsection{Monotonicity of the FE approximation}
\label{sec:mono-fe}
Similar to the continuous problem in \cite{GKSr25}, we prove that 
the discrete QoI \eqref{eq:X_h} is monotone by using the Discrete Maximum Principle (DMP)
to prove a lower bound on the FEM approximation $u_{0, h}$ corresponding to the source term $\ell_0$.

In order to use the DMP we require that the discrete problem
\eqref{eq:fem} is of \emph{nonnegative type} (see, e.g., \cite{CiarRav73}
or Lemma~\ref{lem:nonneg} below), which will be ensured by the following
condition on the sequence of FE meshes. 
First, we introduce some additional notation and parameters of the mesh $\calT_h$.

For a simplex $\tri \in \calT_h$, let $\lambda_k$ for $k = 1, 2, \ldots, d + 1$
denote the barycentric coordinates with respect to the vertices of $\tri$
and define
\[
\sigma(h) \,\coloneqq\,
\max_{\tri \in \calT_h} \sigma(\tri)
\quad \text{with} \quad
\sigma(\tri) \,\coloneqq\, \max_{\substack{k, m = 1, 2, \ldots, d + 1\\ m \neq k}}
\frac{\nabla \lambda_k \cdot \nabla \lambda_m}
{\|\nabla \lambda_k\| \|\nabla \lambda_m\|}\,.
\]

\begin{assumption}
\label{asm:nonneg}
There exists a constant $\sigma_0 < 0$ such that
$\sigma(h) \leq \sigma_0 < 0$ for all $h > 0$ sufficiently small.
\end{assumption}

For example, when $d = 2$ Assumption~\ref{asm:nonneg} holds
if and only if all of the angles in the triangles of the mesh $\calT_h$ 
are less than or equal to $\pi/2$.

\begin{lemma}
\label{lem:nonneg}
Suppose that Assumption~\ref{asm:pde} holds and let $h > 0$ be sufficiently small such that Assumption~\ref{asm:nonneg} holds.
Then for all $\bsy \in \R^{2s + 1}$, the discrete problem \eqref{eq:fem}
is of nonnegative type, i.e., 
\begin{align}
\label{eq:nonneg0}
A_{k, k}(\bsy) \,&>\, 0 \quad \text{for } k = 1, 2, \ldots, \overline{M_h},
\\
\label{eq:nonneg1}
A_{k, m}(\bsy) \,&\leq 0 \quad
\text{for } k = 1, \ldots, M_h, \, m = 1, \ldots, \overline{M_h},\, m \neq k,
\\
\label{eq:nonneg2}
\sum_{m = 1}^{\overline{M_h}} A_{k, m}(\bsy) \,&\geq 0
\quad \text{for } k = 1, 2, \ldots, M_h.
\end{align}
\end{lemma}

\begin{proof}
The proof follows \cite{CiarRav73}, adapting some technical details to fit our setting.

First, by the definition \eqref{eq:stiff} and the bound \eqref{eq:a_bound} we have
\[
A_{k, k}(\bsy) = \int_D a(\bsx, \bsz) |\nabla \phi_k(\bsx)|^2 \, \rd \bsx
\geq \amin(\bsz) \int_D |\nabla \phi_k(\bsx)|^2 \, \rd \bsx
> 0,
\]
proving \eqref{eq:nonneg0}.

To prove \eqref{eq:nonneg1}, we split $A_{k, m}(\bsy)$ into a sum
of integrals on elements $\tri \in \calT_h$ such that
$\tri \subset \supp(\phi_k) \cap \supp(\phi_m) \eqqcolon D_{h, k, m}$,
\begin{align}
\label{eq:nonneg1-sum}
\nonumber
A_{k, m}(\bsy) \,&=\, 
\sum_{\tri \subset D_{h, k, m}} 
\int_\tri a(\bsx, \bsz) \nabla \phi_k(\bsx) \cdot \nabla \phi_m(\bsx) \, \rd \bsx
\\
&=
\sum_{\tri \subset D_{h, k, m}} 
\big(\nabla \phi_k \cdot \nabla \phi_m\big)\big|_\tri \int_\tri a(\bsx, \bsz)  \, \rd \bsx,
\end{align}
where we have use the fact that the gradients $\nabla \phi_k$ and  $\nabla \phi_m$
are piecewise constant on the elements.

Since each $\phi_k$ is piecewise linear, in barycentric coordinates the
restriction to $\tri$ simplifies to $\phi_k|_\tri = \lambda_{k'}$,
where $k' = k'(\tri)$ is the local index on $\tri$ of the $k$th global node.
Hence, the gradients in \eqref{eq:nonneg1-sum} simplify
to give
\[
\big(\nabla \phi_k \cdot \nabla \phi_m\big)\big|_\tri
= \nabla \lambda_{k'} \cdot \nabla \lambda_{m'}
\leq \sigma_0 \|\nabla \lambda_{k'}\| \, \|\nabla \lambda_{m'}\|
\leq \frac{\sigma_0}{\diam(\tri)^2}
\leq \frac{\sigma_0}{h^2} < 0,
\]
where we have used Assumption~\ref{asm:nonneg},
which applies since $k \neq m$ also implies that the local indices
are distinct, $k' \neq m'$.
Substituting this into \eqref{eq:nonneg1-sum} and
using \eqref{eq:a_bound} gives
\[
A_{k, m}(\bsy) \leq
\frac{\sigma_0 \,\amin(\bsz)}{h^2} \sum_{\tri \subset D_{h, k, m}} 
\int_\tri \rd \bsx 
=
\frac{\sigma_0 \, \amin(\bsz)}{h^2} \int_{D_{h, k, m}} \rd \bsx
\leq 0,
\]
as required.

To prove \eqref{eq:nonneg2}, note that the hat functions sum to 1, i.e., 
$\sum_{m = 1}^{\overline{M_h}} \phi_m \equiv 1$ on $\overline{D}$.
Then, by the definition of $A_{k, m}(\bsy)$ in \eqref{eq:stiff},
\[
\sum_{m = 1}^{\overline{M_h}} A_{k, m}(\bsy)
= \sum_{m = 1}^{\overline{M_h}} \calA(\bsy; \phi_k, \phi_m)
= \calA(\bsy; \phi_k, 1) = 0,
\]
since $\nabla 1 \equiv 0$ on $\overline{D}$. 
This implies \eqref{eq:nonneg2}, concluding the proof.
\end{proof}

Since the discrete problem \eqref{eq:fem} is of nonnegative type for all parameters $\bsy \in \R^{2s + 1}$,
the DMP can now be used to show the following lower bound.

\begin{lemma}
\label{lem:u0h_lower}
Suppose that Assumption~\ref{asm:pde} holds, and 
let $h > 0$ be sufficiently small such that Assumption~\ref{asm:nonneg} holds.
Then for all $\bsz \in \R^s$,  the FE solution $u_{0, h}(\cdot, \bsz) \in V_h$ to \eqref{eq:fem}
with right-hand side $\ell_0$ satisfies
\begin{equation}
\label{eq:uh-lower2}
u_{0, h}(\bsx, \bsz) \geq \frac{u_{0, h}(\bsx, \bszero)}{K_{0, h}(\bsz)},
\end{equation}
where, for $K_0(\bsz)$ as in \eqref{eq:K0} and $C_V$ independent of $h$ and $\bsy$,
\begin{equation}
\label{eq:K0h}
K_{0, h}(\bsz) \,\coloneqq\, K_0(\bsz) + \frac{\amax(\bsz) C_Vh^{-d/2}}{\ell_{0, \inf}}.
\end{equation}
\end{lemma}

\begin{proof}
We prove the result by applying the DMP to 
\[
\uhat_h(\bsx, \bsz) \,\coloneqq\, 
\frac{u_{0, h}(\bsx, \bszero)}{K_{0, h}(\bsz)}
 - u_{0, h}(\bsx, \bsz)
 \in V_h,
\]
with coefficient vector $\widehat{\bsU}_h(\bsz) \in \R^{M_h}$. 

First, we show that $\bsA_h(\bsy) \widehat{\bsU}_h(\bsz) \leq 0$, i.e., 
$[\bsA_h(\bsy) \widehat{\bsU}_h(\bsz)]_m \leq 0$ for all $m = 1, 2,\ldots, M_h$.
By the definition of $\bsA_h(\bsy)$ \eqref{eq:stiff} and the expansion \eqref{eq:u_h}
for $\uhat_{h}(\bsz)$ we have
\[
[\bsA_h(\bsy) \widehat{\bsU}_h(\bsz)]_m = \calA(\bsz, \uhat_h(\bsz), \phi_m)
= 
\frac{\calA(\bsz, u_{0, h}(\bszero), \phi_m)}{K_{0, h}(\bsz)}
- \calA(\bsz, u_{0, h}(\bsz), \phi_m) .
\]
The first term can be rewritten and then bounded above using the Cauchy--Schwarz inequality as follows:
\begin{align*}
\calA(\bsz, u_{0, h}(\bszero), \phi_m) \,&=\, 
\calA(\bsz, u_{0}(\bszero), \phi_m) + \calA(\bsz, u_{0, h}(\bszero) - u_0(\bszero), \phi_m)
\\
&\leq K_0(\bsz) \langle \ell_0, \phi_m \rangle
+ \amax(\bsz) \|u_{0, h}(\bszero) - u_0(\bszero)\|_V \, \|\phi_m \|_V,
\end{align*}
where we have also used \cite[eq. (5.5)]{GKSr25} with $v = \phi_m \in V_h \subset V$.
The second term can be bounded using \eqref{eq:u-fe-err-point} (with $\bsz = 0$)
along with the inverse estimate from, e.g., \cite[Theorem 3.2.6]{Ciar02},
and combining these two constants into $C_V$, which is independent of $h$ and $\bsz$,
we obtain
\begin{align}
\label{eq:Au0h<Auh}
\calA(\bsz, u_{0, h}(\bszero), \phi_m) \,&\leq
K_0(\bsz) \langle \ell_0, \phi_m \rangle
+ \amax(\bsz) C_V h \,  h^{-1 - d/2}\int_D |\phi_m(\bsx)|\, \rd \bsx
\nonumber\\
&\leq K_0(\bsz) \langle \ell_0, \phi_m \rangle
+ \amax(\bsz) C_V h^{- d/2} \int_D \frac{\ell_0(\bsx)}{\ell_{0, \inf}} \phi_m(\bsx)\, \rd \bsx
\nonumber\\
&=\, K_{0, h}(\bsz) \langle \ell_0, \phi_m \rangle,
\end{align}
by the definition of $K_{0, h}(\bsz)$ \eqref{eq:K0h}.

Next, using \eqref{eq:fem} and the bound \eqref{eq:Au0h<Auh} gives
\begin{align*}
\calA(\bsz, \uhat_h(\bsz), \phi_m)
\,&=\, 
\frac{\calA(\bsz, u_{0, h}(\bszero), \phi_m)}{K_{0, h}(\bsz)}
 - \calA(\bsz, u_{0, h}(\bsz), \phi_m) 
 \\
&\leq 
\frac{K_{0, h}(\bsz) \langle \ell_0, \phi_m \rangle}{K_{0, h}(\bsz)}
 - \langle \ell_0, \phi_m \rangle 
= 0,
\end{align*}
as required.

Finally, by Lemma~\ref{lem:nonneg} $\bsA_h(\bsz)$ is of nonnegative type
and hence, by the DMP, see e.g.,  \cite[eq.~(2.28)]{CiarRav73}, it follows that
\[
\max_{x \in \overline{D}} \uhat_h(\bsx, \bsz) \,\leq \,
\max_{x \in \partial D} \uhat_h(\bsx, \bsz) = 0,
\]
implying the required result \eqref{eq:uh-lower2}.
\end{proof}

For a positive linear functional $\calG$ this lower bound \eqref{eq:uh-lower2}
implies that the FE approximation of the QoI
is monotone, i.e., 
$\partial^0 \varphi_h(\bsy) = \calG(\partial^0 u_h(\cdot, \bsy)) = \calG(u_{0, h}(\cdot, \bsz)) > 0$
for all $\bsy \in \R^{2s + 1}$. Thus, Assumption~\ref{asm:preint}(a) holds
and the FE preintegrated functions \eqref{eq:g_h-preint} are well defined.

Additionally, Assumption~\ref{asm:preint}(b) also holds by 
taking the limit of \eqref{eq:u_h-linear}.

\subsection{Two lower bounds for the FE approximation}
Lemma~\ref{lem:u0h_lower} from the previous section was useful for showing that the FE approximation is
also monotone with respect to $y_0$. However, the lower bound \eqref{eq:uh-lower2} depends on $h$
and in particular, it tends to 0 as $h \to 0$.
Later, for the FE error analysis we will require explicit lower bounds on 
$\varphi_{0, h}(\bsz) = \calG(u_{0, h}(\bsz))$, and so using the lower bound \eqref{eq:uh-lower2}
would degrade the FE convergence.

Here we present an alternative to Lemma~\ref{lem:u0h_lower}, based on using the lower bound for 
the continuous solution \eqref{eq:phi0-lower} and 
the property that $\varphi_{0, h}(\bsz) \to \varphi_0(\bsz) = \calG(u_0(\cdot, \bsz))$ as $h \to 0$.
It removes the dependence of the lower bound on $h$, but the cost is a ``sufficiently small'' condition
on $h$ that depends on $\bsz$.

\begin{lemma}
\label{lem:lower_fe}
Suppose that Assumption~\ref{asm:pde} holds, let $\bsz \in  \R^{s}$ and
let $\calG \in H^{-1 + \tau}(D)$
for some $\tau \in [0, 1]$ be a positive linear functional.
For $0 < h \leq H(\bsz)$,
the solution $u_{0, h}(\cdot, \bsz) \in V_h$ to the FE problem \eqref{eq:fem} with right hand side
$\ell_0$ satisfies
\begin{equation}
\label{eq:u_h-lower}
\varphi_{0, h}(\bsz) = \calG(u_{0, h}(\cdot, \bsz)) \geq \frac{\calG(u_0(\cdot, \bszero))}{2K_0(\bsz)},
\end{equation}
where
\begin{equation}
\label{eq:H(z)}
H(\bsz) \,\coloneqq\,
\bigg( \frac{\calG(u_0(\cdot, \bszero)) \amin^9(\bsz )}
{2C_\calG \|\ell_0\|_{L^2(D)} K_0(\bsz)   \amax^4(\bsz) \|a(\cdot, \bsz)\|_{C^{0, 1}(\overline{D})}^4}
\bigg)^\frac{1}{1 + \tau}.
\end{equation}
\end{lemma}

\begin{proof}
First, note that since $a(\cdot, \bsz) \in W^{1, \infty}(D)$ and $D$ is bounded it follows that
$a(\cdot, \bsz) \in C^{0, 1}(\overline{D})$
and clearly $\|a(\cdot, \bsz)\|_{C^{0, \tau}(\overline{D})} \leq \|a(\cdot, \bsz)\|_{C^{0, 1}(\overline{D})}$. Furthermore, since $\ell_0 \in L^\infty(D)  \subset L^2(D)$ (because $D$ is bounded) whereas
$\calG \in H^{-1 + \tau}(D)$,
the FE error bound \eqref{eq:fe-err-qoi-point}
becomes
\begin{equation}
\label{eq:fe-err-L2}
\big|\calG(u_0(\cdot, \bsz)) - \calG(u_{0,h}(\cdot, \bsz))\big|
\leq
C_\calG \frac{\amax^4(\bsz) \|a(\cdot, \bsz)\|_{C^{0, 1}(\overline{D})}^4}{\amin^9(\bsz)}
\|\ell_0\|_{L^2(D)}
\, h^{1 + \tau}\,.
\end{equation}
From which it follows that for all $h \leq H(\bsz)$, with $H(\bsz)$ as in \eqref{eq:H(z)},
we have
\begin{equation}
\label{eq:qoi-H-bound}
\big|\calG(u_0(\cdot, \bsz)) - \calG(u_{0, h}(\cdot, \bsz))\big| \leq
\frac{\calG(u_0(\cdot, \bszero))}{2K_0(\bsz)}\,.
\end{equation}

If $\calG(u_{0, h}(\cdot, \bsz)) \geq \calG(u_0(\cdot,\bsz)) > 0$
then the result holds trivially
by the linearity of $\calG$ and the lower bound \eqref{eq:phi0-lower} on
$\calG(u(\cdot, \bsz))$.
For the remainder of the proof suppose $\calG(u_{0, h}(\cdot, \bsz)) \leq \calG(u_0(\cdot, \bsz))$.
Under this assumption, we have
\begin{align*}
\calG(u_{0, h}(\cdot, \bsz)) \,&=\,
\calG(u_0(\cdot,\bsz)) - \big[\calG(u_0(\cdot, \bsz)) - \calG(u_{0,h}(\cdot, \bsz))\big]
\\
&=\, \calG(u_0(\cdot, \bsz)) - \big|\calG(u_0(\cdot, \bsz)) - \calG(u_{0, h}(\cdot, \bsz))\big|
\\
&\geq \frac{ \calG(u_0(\cdot, \bszero))}{K_0(\bsz)} 
- \frac{ \calG(u_0(\cdot, \bszero))}{2K_0(\bsz)}
=\frac{ \calG(u_0(\cdot, \bszero))}{2K_0(\bsz)},
\end{align*}
where in the last inequality we have substituted in \eqref{eq:phi0-lower}
and \eqref{eq:qoi-H-bound}.
\end{proof}

Note that the restriction on the meshwidth in Lemma~\ref{lem:lower_fe}, i.e., $h \leq H(\bsz)$, depends on $\bsz$ and
in particular, $H(\bsz) \to 0$ if any $z_j \to \pm \infty$. This comes from the
fact that as $z_j \to \pm \infty$ the problem becomes more ill conditioned,
with the lower bound on the continuous solution \eqref{eq:phi0-lower} tending to 0
and the FE approximation becoming less accurate, since the constant in \eqref{eq:fe-err-qoi-point}
tends to $\infty$. Hence, there is no value of $h > 0$
such that the required lower bound holds for all $\bsz \in \R^{s}$.

A second lower bound can be obtained from Lemma~\ref{lem:u0h_lower},
which removes this condition on $h$, but at the cost of a factor $h^{d/2}$ in the lower bound,
implying that the lower bound tends to 0 as $h \to 0$.

\begin{corollary}
Suppose that Assumption~\ref{asm:pde} holds,
let $h > 0$ be sufficiently small such that Assumption~\ref{asm:nonneg} holds
and let $\calG \in V^*$ be a  positive linear functional. 
Then, for all $\bsz \in \R^s$, the FE solution $u_{0, h}(\cdot, \bsz) \in V_h$
to \eqref{eq:fem} with right-hand side $\ell_0$ satisfies
\begin{equation}
\label{eq:phi_h-lower2}
\varphi_{0, h}(\bsz) = \calG(u_{0, h}(\cdot, \bsz)) 
\geq \frac{\amin(\bszero)}{\amax(\bszero)} \frac{\calG(u_{0}(\cdot, \bszero))}{2K_{0, h}(\bsz)} \simeq h^{d/2},
\end{equation}
with $K_{0, h}(\bsz)$ as in \eqref{eq:K0h}.
\end{corollary}

\begin{proof}
Since $\calG$ is positive, it follows from \eqref{eq:uh-lower2} that
\[
\calG(u_{0, h}(\cdot, \bsz)) 
\geq \frac{\calG(u_{0, h}(\cdot, \bszero))}{K_{0, h}(\bsz)}\,.
\]
For $h$ sufficiently small such that $h \leq H(\bszero)$, the lower bound \eqref{eq:u_h-lower} for $\bsz = \bszero$
becomes
\[
\calG(u_{0, h}(\cdot, \bszero)) \geq \frac{\calG(u_{0}(\cdot, \bszero))}{2C_0(\bszero)} = 
\frac{\amin(\bszero)}{\amax(\bszero)} \frac{\calG(u_{0}(\cdot, \bszero))}{2}
\]
where we have used that $K_0(\bszero) = \amax(\bszero)/\amin(\bszero)$.
The result \eqref{eq:phi_h-lower2} follows by combining the two bounds above.
\end{proof}

In summary, here we have shown two lower bounds on $\varphi_{0, h}(\bsz)$.
The first bound \eqref{eq:u_h-lower} is similar to the lower bound \eqref{eq:phi0-lower} for the continuous solution,
but needs $h$ sufficiently small depending on $\bsz$,
whereas the second bound \eqref{eq:phi_h-lower2} holds for all $\bsz$ and $h$ sufficiently small 
(independent of $\bsz$), but depends poorly on $h$.
Both bounds will be useful for the FE error analysis in Section~\ref{sec:fe-err}.

\subsection{Parametric regularity of the FE approximation}

Since $V_h$ is conforming, i.e., $V_h \subset V$, the FE approximations have the same parametric regularity as the continuous versions.
In particular, the derivative bounds from \cite[Theorem 3.1]{GKSr25} also hold for $\|\partial^\bsnu u_h(\cdot, \bsy)\|_V$
and thus Assumption~\ref{asm:preint}(c) holds for any $\bsnu \in \N_0^{2s}$.
A more detailed analysis of the parametric regularity of the FEM approximation $u_h(\cdot, \bsy)$,
including explicitly bounding the $\calH_{2s}^\bsnu$-norm of the preintegrated functions in \eqref{eq:g_h-preint},
can also be performed by following Section 5.3 of \cite{GKSr25}.
Indeed, \cite[Lemma~5.6]{GKSr25} also holds for the FE counterparts by replacing $\varphi_0(\bsz)$ by
$\varphi_{0, h}(\bsz)$ in the upper bound and then a similar bound as in \cite[Lemma 5.7]{GKSr25} also
holds by instead using one of the lower bounds on $\varphi_{0, h}(\bsz)$ from Section~\ref{sec:mono-fe},
which in turn can be used to bound the norm as in \cite[eq. (8.1)]{GKSr25}.
However, such detail is not required for the error analysis to follow.

\section{Error analysis}
\label{sec:error}

\subsection{FE error}
\label{sec:fe-err}
Now we bound the error of the FE approximations of the cdf \eqref{eq:cdf0} and pdf \eqref{eq:pdf0}.
The analysis relies on smoothing effect of preintegration to remove the discontinuities
in the original formulations, making use of the smooth preintegrated formulations \eqref{eq:F_h-preint}
and \eqref{eq:f_h-preint} with \eqref{eq:g_h-preint}. First, we require the following technical results.

\begin{lemma}
The standard normal cdf $\Phi$ and pdf $\rho$ satisfy,
for $y, y' \in \R$,
\begin{align}
\label{eq:Lip-Phi}
|\Phi(y) - \Phi(y')| \,&\leq \frac{1}{\sqrt{2\pi}} |y - y'|,
\qquad \text{and}
\\
\label{eq:rho-diff}
|\rho(y) - \rho(y')| \,&\leq \frac{1}{\sqrt{2\pi}} \max(|y|, |y'|)  |y - y'|.
\end{align}
\end{lemma}
\begin{proof}
Without loss of generality assume $y > y'$. By the Mean Value Theorem (MVT),
for some $\chi \in (y', y)$,
\[
\Phi(y) - \Phi(y') = \rho(\chi)(y - y'),
\]
where we have used that $\Phi' = \rho$. Taking the absolute value
and then using the property that
$\max_{\chi \in \R} \rho(\chi) = \rho(0) = 1/\sqrt{2\pi}$
gives the result \eqref{eq:Lip-Phi}.

Again using the MVT on $\rho$ now, for some $\chi \in (y', y)$,
\[
\rho(y) - \rho(y') = \rho'(\chi)(y - y') =
-\chi \rho(\chi) (t - t').
\]
Taking the absolute value gives
\[
|\rho(y) - \rho(y')| =
|\chi| \rho(\chi) |y - y'|
\leq \frac{1}{\sqrt{2\pi}} \max( |y|, |y'|) |y - y'|,
\]
where we have used again the upper bound on $\rho$ and the fact that
$\chi \in (y', y)$
so that $|\chi| \leq \max (|y|, |y'|)$.
\end{proof}
\begin{lemma}
Let Assumption~\ref{asm:pde} hold, then for all $\bsx \in D$ and  $\bsz \in \R^s$
the coefficient $a(\bsx, \bsz)$ satisfies the following bounds
\begin{align}
\label{eq:amax-bnd}
a(\bsx, \bsz) \,&\leq \amax(\bsz) \leq
\exp\big(\|a_0\|_{L^\infty(D)}\big)\exp\Bigg( \sum_{j = 1}^s |z_j| \|a_j\|_{L^\infty(D)}\Bigg),
\\
\label{eq:amin-bnd}
a(\bsx, \bsz) \,&\geq \amin(\bsz) \geq \exp\Bigg(\inf_{\bsx \in D} a_0(\bsx)\Bigg)
\exp\Bigg( -\sum_{j = 1}^s |z_j| \|a_j\|_{L^\infty(D)}\Bigg),
\end{align}
and
\begin{equation}
\label{eq:||a||-bnd}
\|a(\cdot, \bsz) \|_{C^{0, 1}(\overline{D})} \leq
\amax(\bsz)  \bigg( \|a_0\|_{C^{0, 1}(\overline{D})} +
\sum_{j = 1}^s |z_j| \|a_j\|_{C^{0, 1}(\overline{D})}\bigg).
\end{equation}
\end{lemma}
\begin{proof}
The first two bounds follow easily from the definition of $a$ \eqref{eq:coeff}.
Next we bound $\|a(\cdot, \bsz)\|_{C^{0, 1}(\overline{D})}$. Let $\bsx, \bsx' \in D$,
by applying the mean value theorem to $\exp(\cdot)$, for some $\chi$ between
$a_0(\bsx) + \sum_{j = 1}^s z_j  a_j(\bsx)$ and $a_0(\bsx') + \sum_{j = 1}^s z_j  a_j(\bsx')$
we have
\begin{align*}
|a(\bsx, \bsy) - a(\bsx', \bsy) |
\,&=\, \exp(\chi)\bigg| a_0(\bsx) + \sum_{j = 1}^s z_j  a_j(\bsx) - a_0(\bsx') - \sum_{j = 1}^s z_j  a_j(\bsx')\bigg|
\\
&\leq \amax(\bsz) \bigg| a_0(\bsx)  - a_0(\bsx')
+ \sum_{j = 1}^s z_j  (a_j(\bsx) - a_j(\bsx'))\bigg|\,.
\end{align*}
Dividing through by $|\bsx - \bsx'|$, using the triangle inequality
 and then taking the $\sup$ over $\bsx, \bsx' \in D$ we obtain the bound \eqref{eq:||a||-bnd}.
\end{proof}

The two lower bounds on the FE solution \eqref{eq:u_h-lower} and \eqref{eq:phi_h-lower2}
will play an important role in the FE error analysis of the cdf and pdf. 
To utilise both bounds, we will split $\R^s$ into a bounded box and its complement,
then use a different lower bound on each subdomain. 
The following lemma shows that the error in truncating the domain to this box is
of the correct order for our FE analysis.
We use the notation $\abs(\bsz)$ below for the componentwise absolute
value of a vector. 

\begin{lemma}
\label{lem:Omega_h}
Suppose that Assumption~\ref{asm:pde} holds, let $h$ be sufficiently 
small and let $\tau$, $\tau' \in [0, 1]$.
For $\bsbeta \coloneqq (\beta_j)_{j = 1}^s \in [0, \infty)^s$, 
define the box $\Omega_h = \Omega_h(\bsbeta) \subset \R^s$ by
\begin{equation}
\label{eq:Omega_h}
\Omega_h(\bsbeta) \,\coloneqq\, \bigotimes_{j = 1}^s [-\kappa_h - \beta_j, \kappa_h + \beta_j] 
\end{equation}
where, for $c_{\bsbeta}$ independent of $h$ and $\bsz$ defined in \eqref{eq:c_beta},
\begin{equation}
\label{eq:kappa}
\kappa_h \,\coloneqq\, -\Phi^{-1} \bigg(\frac{(c_{\bsbeta} h^{d/2 + \tau + \tau'})^{\frac{1}{s}}}{2}\bigg)
= \bigg|\Phi^{-1} \bigg(\frac{(c_{\bsbeta} h^{d/2 + \tau + \tau'})^{\frac{1}{s}}}{2}\bigg)\bigg|\,.
\end{equation}
Then $h \leq H(\bsz)$ for all $\bsz \in \Omega_h(\bsbeta)$ and, for $C_{\bsbeta, \tau} < \infty$ independent of $h$ and $\bsz$,
\begin{equation}
\label{eq:int_not_Omega_h}
\int_{\R^s \setminus \Omega_h }
\exp\big(\abs(\bsz) \cdot \bsbeta\big) \bsrho(\bsz) \, \rd \bsz
\leq C_{\bsbeta, \tau} h^{d/2 + \tau + \tau'}\,.
\end{equation}
\end{lemma}

\begin{proof}
We first show that $h \leq H(\bsz)$, then we prove the bound on the integral \eqref{eq:int_not_Omega_h}.
Assume that $h$ is sufficiently small such that $\kappa_h \geq 1$.

Recall from \eqref{eq:H(z)} that
\begin{align}
\label{eq:H0}
&H(\bsz)^{1 + \tau} =
\frac{\ellinf \,\calG(u_0(\cdot, \bszero)) }{2C_\calG \|\ell_0\|_{L^2(D)} }
\\\nonumber
&\cdot 
\underbrace{\frac{\amin^9(\bsz)}
{\amax^5(\bsz) \|a(\bsz)\|_{C^{0, 1}(\overline{D})}^4
\big[ \ellinf + \|u_0(\cdot, \bszero)\|_{W^{1, \infty}}
\sum_{j = 1}^s |z_j|\|a_j\|_{W^{1, \infty}(D)}\big] }}_{\widetilde{H}(\bsz)},
\end{align}
where we have expanded $K_0(\bsz)$ and
simplified into two factors, the first of which is independent of $\bsz$
and the second we denote by $\widetilde{H}(\bsz)$.

By first deriving a lower bound on $\widetilde{H}(\bsz)$, we now show that the box
$\Omega_h$ is such that
$h \leq H(\bsz)$ for all $\bsz \in \Omega_h$.
From the three bounds \eqref{eq:amax-bnd}--\eqref{eq:||a||-bnd} above it follows that $\widetilde{H}(\bsz)$ is bounded from below by
\begin{align*}
\widetilde{H}(\bsz) \geq
&\exp\bigg(9\inf_{\bsx \in D} a_0(\bsx) - 9\|a_0\|_{L^\infty} \bigg)
\cdot
\frac{\exp\big(-18\sum_{j = 1}^s |z_j| \|a_j\|_{L^\infty}\big)}
{\big[\|a_0\|_{C^{0, 1}(\overline{D})} + \sum_{j = 1}^s |z_j| \|a_j\|_{C^{0, 1}(\overline{D})}\big]^4 }
\\
&\cdot
\frac{1}
{ \ellinf + \|u_0(\cdot, \bszero)\|_{W^{1, \infty}}
\sum_{j = 1}^s |z_j|\|a_j\|_{W^{1, \infty}} }\,.
\end{align*}
For $\bsz \in \Omega_h$, this lower bound is clearly minimised on the boundary,
i.e., when $z_j = \pm (\kappa_h + \beta_j)$ for $j = 1, 2, \ldots, s$.
In this case,
\begin{align*}
\widetilde{H}(\bsz)
\geq
&\exp\bigg(9\inf_{\bsx \in D} a_0(\bsx) - 9\|a_0\|_{L^\infty} \bigg)
\cdot
\frac{\exp\big(-18 \sum_{j = 1}^s (\kappa_h + \beta_j)\|a_j\|_{L^\infty}\big)}
{\big[\|a_0\|_{C^{0, 1}(\overline{D})} + \sum_{j = 1}^s (\kappa_h + \beta_j)\|a_j\|_{C^{0, 1}(\overline{D})}\big]^4 }
\\
&\cdot
\frac{1}
{ \ellinf + \|u_0(\cdot, \bszero)\|_{W^{1, \infty}}
\sum_{j = 1}^s (\kappa_h + \beta_j)\|a_j\|_{W^{1, \infty}} }
\\
\geq
&\frac{\exp\big(-18 \kappa_h \sum_{j = 1}^s\|a_j\|_{L^\infty}\big)}{\kappa_h^5} \\
&\cdot\frac{ \exp \big(9\inf_{\bsx \in D} a_0(\bsx) - 9\|a_0\|_{L^\infty}
-18 \sum_{j = 1}^s\beta_j \|a_j\|_{L^\infty}\big)}
{\big[\|a_0\|_{C^{0, 1}(\overline{D})} + \sum_{j = 1}^s (1 + \beta_j)\|a_j\|_{C^{0, 1}(\overline{D})}\big]^4 }
\\
&\cdot
\frac{1}
{ \ellinf + \|u_0(\cdot, \bszero)\|_{W^{1, \infty}}
\sum_{j = 1}^s (1 + \beta_j)\|a_j\|_{W^{1, \infty}} }\,.
\end{align*}
where we have used the fact that $ \kappa_h \geq 1$.

Using the asymptotic formula $-\Phi^{-1}(\chi) = \Phi^{-1}(1 - \chi) \sim \sqrt{2\log(1/\chi)}$ from \cite{BEJ76}, it follows that
$\Phi^{-1}(\chi) \leq 2 \sqrt{\log(1/\chi)}$ for $\chi$ sufficiently small.
Thus, for $h$ sufficiently small (independent of $\bsy$), letting $\alpha_0 = 18 \sum_{j = 1}^s \|a_j\|_{L^\infty(D)}$
and $\chi = (c_{\bsbeta} h^{d/2 + \tau + \tau'})^{\frac{1}{s}}/2$, the first factor above 
can be bounded from below by
\[
\frac{\exp( -\alpha_0 \kappa_h)}{\kappa_h^5}
\geq \frac{\exp( -\sqrt{\log(1/\chi)})^{2\alpha_0}}{2^5\log(1/\chi)^{5/2}}
\geq \frac{\chi^{2\alpha_0/\varepsilon}}{2^5 \log(1/\chi)^{5/2}}
\geq \frac{\chi^{2\alpha/\varepsilon}}{2^5},
\]
where we have also used that for any $\varepsilon > 0$ we have $\exp(-\sqrt{\log(1/\chi)}) \geq \chi^\varepsilon$
for $\chi$ sufficiently small,
which follows from the property that for $k = 1/\chi$ sufficiently large
$\exp(\sqrt{\log(k)}) \leq k^{\varepsilon}$ for $\varepsilon > 0$.
Then in the last step we take $\alpha > \alpha_0 = 18 \sum_{j = 1}^s \|a_j\|_{L^\infty(D)}$
to absorb the log factor in the denominator. 

Choosing now $\varepsilon = 2\alpha  (d/2 + \tau + \tau') / [s(1 + \tau)]  > 0$ gives
\[
\frac{\exp\big(-18 \kappa_h \sum_{j = 1}^s\|a_j\|_{L^\infty}\big)}{\kappa_h^5}
\geq 
\frac{\chi^{s(1 + \tau)/(d/2 + \tau + \tau')}}{2^5}
= \frac{c_{\bsbeta}^{(1 + \tau)/(d/2 + \tau + \tau')} h^{1 + \tau}}{2^{5 + s(1 + \tau)/(d/2 + \tau + \tau')}}.
\]
Thus, for $\bsz \in \Omega_h$, we have
\begin{align*}
\widetilde{H}(\bsz)
&\geq  \frac{c_{\bsbeta}^{(1 + \tau)/(d/2 + \tau + \tau')} h^{1 + \tau}}{2^{5 + s(1 + \tau)/(d/2 + \tau + \tau')}}
\cdot \frac{9 \sum_{j = 1}^s \|a_j\|_{L^\infty}}
{\big[\|a_0\|_{C^{0, 1}(\overline{D})} + \sum_{j = 1}^\infty (1 + \beta_j) \|a_j\|_{C^{0, 1}(\overline{D})}\big]^4}
\\
&\cdot
\frac{\exp\big(9\inf_{\bsx \in D} a_0(\bsx) - 9\|a_0\|_{L^\infty(D)} \big)}
{ \ellinf + \|u_0(\cdot, \bszero)\|_{W^{1, \infty}}
\sum_{j = 1}^s (1 + \beta_j)\|a_j\|_{W^{1, \infty}} }\,.
\end{align*}

Substituting this into the formula \eqref{eq:H0} for $H(\bsz)$ gives
\begin{align*}
H(\bsz) &\geq
\frac{c_{\bsbeta}^{1/(d/2 + \tau + \tau')} h}{2^{s/(d/2 + \tau + \tau')}}
\Bigg( \frac{9 \,\ellinf \,\calG(u_0(\cdot, \bszero))
\sum_{j = 1}^s \|a_j\|_{L^\infty}}
{2^6 C_\calG \|\ell_0\|_{L^2(D)} \big[\|a_0\|_{C^{0, 1}(\overline{D})} + \sum_{j = 1}^s (1 + \beta_j)\|a_j\|_{C^{0, 1}(\overline{D})}\big]^4}
\\
&\qquad\cdot
\frac{\exp\big(9\inf_{\bsx \in D} a_0(\bsx) - 9\|a_0\|_{L^\infty(D)} \big)}
{\ellinf + \|u_0(\cdot, \bszero)\|_{W^{1, \infty}}
\sum_{j = 1}^s (1 + \beta_j)\|a_j\|_{W^{1, \infty}}}
\Bigg)^\frac{1}{1 + \tau}
\\
&=  h, 
\end{align*}
where we now define $c_{\bsbeta}$ by
\begin{align}
\label{eq:c_beta}
c_{\bsbeta} &\coloneqq 2^s
\Bigg( \frac{9 \,\ellinf \,\calG(u_0(\cdot, \bszero))
\sum_{j = 1}^s \|a_j\|_{L^\infty}}
{2^6C_\calG \|\ell_0\|_{L^2(D)} \big[\|a_0\|_{C^{0, 1}(\overline{D})} + \sum_{j = 1}^s (1 + \beta_j)\|a_j\|_{C^{0, 1}(\overline{D})}\big]^4}
\nonumber\\
&\cdot
\frac{\exp\big(9\inf_{\bsx \in D} a_0(\bsx) - 9\|a_0\|_{L^\infty(D)} \big)}
{\ellinf + \|u_0(\cdot, \bszero)\|_{W^{1, \infty}}
\sum_{j = 1}^s (1 + \beta_j)\|a_j\|_{W^{1, \infty}}}
\Bigg)^{(d/2 + \tau + \tau')/(1 + \tau)},
\end{align}
which is clearly independent of $h$ and $\bsz$.
Hence, $h \leq H(\bsz)$ for all $\bsz \in \Omega_h$.

Finally, we prove the bound \eqref{eq:int_not_Omega_h}
\begin{align}
&\int_{\R^{s} \setminus \Omega_h}
\exp\Bigg(\sum_{j = 1}^s |z_j| \beta_j\Bigg) \bsrho(\bsz) \, \rd \bsz
\nonumber\\
&= \prod_{j = 1}^{s} \bigg[\int_{-\infty}^{-\kappa_h - \beta_j} 
\frac{\exp\big(-z_j^2/2 + \beta_j |z_j|\big)}{\sqrt{2\pi}}\, \rd z_j
+ \int_{\kappa_h + \beta_j}^\infty \frac{\exp\big(-z_j^2/2 + \beta_j |z_j|\big)}{\sqrt{2\pi}}\, \rd z_j\bigg]
\nonumber\\
&= \prod_{j = 1}^{s} 2\int_{-\infty}^{-\kappa_h - \beta_j} 
\frac{\exp\big(-z_j^2/2 - \beta_j z_j\big)}{\sqrt{2\pi}}\, \rd z_j
\nonumber\\
&= \prod_{j = 1}^{s} 2\exp(\beta_j^2/2) \int_{-\infty}^{-\kappa_h - \beta_j} 
\rho(z_j + \beta_j)\, \rd z_j
\nonumber\\
&= \exp\Bigg(\sum_{j = 1}^s \frac{\beta_j^2}{2}\Bigg)\prod_{j = 1}^{s} 2 \int_{-\infty}^{-\kappa_h - \beta_j} 
\rho(z_j + \beta_j)\, \rd z_j
\nonumber\\
&=  \exp(\bsbeta \cdot \bsbeta/2)
\big[ 2\Phi(-\kappa_h)\big]^s 
= c_{\bsbeta}  \exp(\bsbeta \cdot \bsbeta/2) h^{d/2 + \tau + \tau'} ,
\end{align}
which gives the desired result with $C_{\bsbeta, \tau} = c_{\bsbeta}  \exp(\bsbeta \cdot \bsbeta/2)$ independent of $h$ and $\bsz$.
\end{proof}

Now we are ready to bound the FE error of the cdf and pdf approximations.
The strategy is to use the simplifications of $F_h$ and $f_h$
after performing preintegration as in \eqref{eq:F_h-preint} and \eqref{eq:f_h-preint}, respectively.

\begin{theorem}
\label{thm:fe-err}
Suppose that Assumption~\ref{asm:pde} holds, for all $\bsw \in \R^{s + 1}$
let $\ell(\cdot, \bsw) \in H^{-1 + \tau'}(D)$ for some $\tau' \in (0, 1]$
and let $\calG \in H^{-1 + \tau}(D)$ for some $\tau \in [0, 1]$ be a positive linear functional.
For $h > 0$ sufficiently small and $t \in [t_0, t_1]$, the FE approximations $F_{h}(t)$ \eqref{eq:F_h-preint}
and $f_{h}(t)$ \eqref{eq:f_h-preint} of the cdf
and pdf of the quantity of interest $X = \calG(u(\cdot, \bsy))$ satisfy
\begin{align}
\label{eq:F-fe-err}
|F(t) - F_{h}(t)|
\,&\leq C_{F, t, \tau} \, h^{\tau' + \tau},
\quad \text{and}\\
\label{eq:f-fe-err}
|f(t) - f_{h}(t)|
\,&\leq C_{f, t, \tau} \, h^{\tau' + \tau},
\end{align}
for $C_{F, t, \tau}, C_{f, t, \tau} < \infty$ independent of $h$ and $\bsy$.
\end{theorem}

\begin{proof}
Using \eqref{eq:F_h-preint} with \eqref{eq:g_h-preint} and the corresponding simplification for $F(t)$,
we can write the FE error for the cdf as
\begin{align}
|F(t) - F_h(t)| &= \bigg|\int_{\R^{2s}} \big[\Phi(\xi(t, \bsy)) - \Phi(\xi_h(t, \bsy))\big]
\bsrho(\bsy) \, \rd \bsy\bigg| 
\nonumber\\
\label{eq:F_h-err}
&\leq
\int_{\R^{2s}} \big|\Phi(\xi(t, \bsy)) - \Phi(\xi_h(t, \bsy))\big|
\bsrho(\bsy) \, \rd \bsy.
\end{align}

Next, we separate the variable $\bsy$ into $\bsw \in \R^s$ and $\bsz \in \R^s$,
then split the domain of integration for $\bsz \in \R^s$ 
into a bounded box $\Omega_F \coloneqq \Omega_h(\bszero)$
(cf.~\eqref{eq:Omega_h} with $\bsbeta = \bszero$) and its exterior to give
\begin{align}
\label{eq:cdf_int_split}
|F(t) - F_h(t)|
\leq&
\underbrace{\int_{\R^s}\int_{\Omega_F} \big|\Phi(\xi(t, \bsw, \bsz))
- \Phi(\xi_h(t, \bsw, \bsz))\big|
\bsrho(\bsw)\bsrho(\bsz) \rd \bsw\rd \bsz
}_{I_{F, 1}}
\nonumber\\&
+ \underbrace{\int_{\R^s}\int_{\R^s \setminus \Omega_F} \big|\Phi(\xi(t, \bsw, \bsz))
- \Phi(\xi_h(t, \bsw, \bsz))\big|
\bsrho(\bsw)\bsrho(\bsz) \, \rd \bsw\rd\bsz
}_{I_{F, 2}}.
\end{align}

The first integral can be simplified using the Lipschitz bound \eqref{eq:Lip-Phi} as follows
\begin{equation}
\label{eq:F_h-err-Lip}
I_{F, 1}
\leq  \frac{1}{\sqrt{2\pi}} \int_{\R^s}\int_{\Omega_F}
|\xi(t, \bsw, \bsz) - \xi_h(t, \bsw, \bsz) | \bsrho(\bsw)\bsrho(\bsz) \,\rd\bsw \rd\bsz.
\end{equation}
Expanding $\xi_h$ using \eqref{eq:xi_h} along with the corresponding formula for $\xi$, 
then rearranging and using the triangle inequality we obtain the bound
\begin{align*}
|\xi(t, \bsy) - \xi_h(t, \bsy) |
\leq
&\frac{|\phibar(\bsz) - \phibar_h(\bsz)| +
\sum_{i = 1}^s |w_i|\, |\varphi_i(\bsy) - \varphi_{i, h}(\bsy)|}{\varphi_0(\bsz)}
\\
&+ \frac{|\varphi_{0}(\bsy) - \varphi_{0, h}(\bsy)| \, |t - \phibar_h(\bsz) - \sum_{i = 1}^s w_i \varphi_{i, h}(\bsz)|}
{\varphi_0(\bsz) \varphi_{0, h}(\bsz)}\,.
\end{align*}
For $\bsy = (\bsw, \bsz) \in \R^{s} \times \Omega_F$, 
substituting in the FE error bound \eqref{eq:fe-err-qoi-point} corresponding to
the right-hand sides $\ellbar \in H^{-1 + \tau'}(D)$ and $\ell_i \in H^{-1 + \tau'}(D)$,
the lower bounds \eqref{eq:phi0-lower}
and \eqref{eq:u_h-lower},
as well as the upper bound
$\phibar_{h}(\bsz) \leq \|\calG\|_{V^*} \|\ellbar\|_{V^*} / \amin(\bsz)$
(which follows from \eqref{eq:LaxMil-fe} with right-hand side $\ellbar$) and the corresponding bounds
for $\varphi_{i, h}(\bsz)$,
we have
\begin{align}
\label{eq:xi-fe-err}
|\xi(t, \bsw, \bsz) &- \xi_h(t, \bsw, \bsz) |
\nonumber\\
\leq&
\frac{C_\calG\, K_0(\bsz)\, \amax^4(\bsz) \|a(\cdot, \bsz)\|_{C^{0, 1}(\overline{D})}^4}
{\varphi_0(\bszero)  \amin^9(\bsz)}
\bigg( \|\ellbar\|_{L^2(D)} + \sum_{i = 1}^s |w_i| \|\ell_i\|_{L^2(D)}
\nonumber\\&
+ \frac{2\|\ell_0\|_{L^2(D)} K_0(\bsz)}{\varphi_0(\bszero)}
\bigg[ |t| +  \frac{\|\calG\|_{V^*}\big[ \|\ellbar\|_{V^*} + \sum_{i = 1}^s |w_i| \|\ell_i\|_{V^*}\big]}{\amin(\bsz)}\bigg]
\bigg)\,
h^{\tau' + \tau}\,.
\end{align}
Note that Lemma~\ref{lem:Omega_h} ensures that $h \leq H(\bsz)$
for all $\bsz \in \Omega_F$,
implying that the lower bound \eqref{eq:u_h-lower} holds
for all $\bsz \in \Omega_F$.

Substituting the bound \eqref{eq:xi-fe-err} into \eqref{eq:F_h-err-Lip} gives 
\begin{align*}
I_{F, 1}
\leq&
\frac{C_\calG h^{\tau' + \tau}}{\sqrt{2\pi} \varphi_0(\bszero)}
\int_{\R^s} \!\int_{\Omega_F}  
\frac{ K_0(\bsz) \amax^4(\bsz) \|a(\cdot, \bsz)\|_{C^{0, 1}(\overline{D})}^4}
{  \amin^9(\bsz)}
\bigg(\! \|\ellbar\|_{L^2(D)} + \sum_{i = 1}^s |w_i| \|\ell_i\|_{L^2(D)}
\\&
+ \frac{2\|\ell_0\|_{L^2(D)} K_0(\bsz)}{\varphi_0(\bszero)}
\bigg[ |t| +  \frac{\|\calG\|_{V^*}\big[ \|\ellbar\|_{V^*} + \sum_{i = 1}^s |w_i| \|\ell_i\|_{V^*}\big]}{\amin(\bsz)}\bigg]
\bigg) \bsrho(\bsw)\bsrho(\bsz) \, \rd \bsw \rd \bsz
\\
\leq&
\frac{C_\calG h^{\tau' + \tau}}{\sqrt{2\pi} \varphi_0(\bszero)}
\int_{\R^{s}}\int_{\R^s} \frac{ K_0(\bsz)\amax^4(\bsz) \|a(\cdot, \bsz)\|_{C^{0, 1}(\overline{D})}^4}
{  \amin^9(\bsz)}
\bigg( \|\ellbar\|_{L^2(D)} + \sum_{i = 1}^s |w_i| \|\ell_i\|_{L^2(D)}
\\&
+ \frac{2\|\ell_0\|_{L^2(D)} K_0(\bsz)}{\varphi_0(\bszero)}
\bigg[ |t| +  \frac{\|\calG\|_{V^*}\big[ \|\ellbar\|_{V^*} + \sum_{i = 1}^s |w_i| \|\ell_i\|_{V^*}\big]}{\amin(\bsz)}\bigg]
\bigg) \bsrho(\bsw)\bsrho(\bsz) \, \rd \bsw \rd \bsz
\\
=& C_{F, 1} \, h^{\tau' + \tau},
\end{align*}
where the integral is finite by using the bounds \eqref{eq:amax-bnd}--\eqref{eq:||a||-bnd}.

Next, using the fact the $\Phi$ is bounded above by 1
and then using the upper bound \eqref{eq:int_not_Omega_h} (with $\bsbeta = \bszero$),
we can bound the second integral in \eqref{eq:cdf_int_split} by
\begin{align*}
I_{F, 2} \,&\leq
2\int_{\R^s}\int_{ \R^s \setminus \Omega_F} \bsrho(\bsw)\bsrho(\bsz) \,\rd\bsw\rd \bsz
= 2\int_{\R^s} \bsrho(\bsw)\,\rd\bsw
\int_{\R^s\setminus \Omega_h(\bszero)} \bsrho(\bsz)\,\rd\bsz
\\
&\leq 2C_{\bszero, \tau} h^{d/2 + \tau + \tau'}
\leq C_{F, 2} h^{\tau + \tau'},
\end{align*}
where we have used the fact that $h$ is sufficiently small.
Substituting the bounds on $I_{F, 1}$ and $I_{F, 2}$
into \eqref{eq:cdf_int_split} gives the required result \eqref{eq:F-fe-err}.

Similarly, using \eqref{eq:f_h-preint} with \eqref{eq:g_h-preint} and the corresponding simplification for $f(t)$,
we can bound the FE error for the pdf by
\begin{equation}
\label{eq:f_h-err1}
|f(t) - f_h(t)| \leq \int_{\R^{2s}} \bigg|
\frac{\rho(\xi(t, \bsy))}{\varphi_0(\bsz)}
 -
\frac{\rho(\xi_h(t, \bsy))}{\varphi_{0, h}(\bsz)}
\bigg| \bsrho(\bsy) \, \rd \bsy\,.
\end{equation}
Again, separating $\bsy = (\bsw, \bsz)$ and splitting the domain of integration for $\bsz \in \R^s$
based on the box $\Omega_f \coloneqq \Omega_h(\bsbeta)$ with $\bsbeta$ given by $\beta_j  \coloneqq \|a_j\|_{L^\infty} + \|a_j\|_{W^{1, \infty}}$, 
we can write this as
\begin{align}
\label{eq:f_h-err2}
|f(t) - f_h(t)|
&\leq
\underbrace{\int_{\R^s}\int_{\Omega_f} \bigg|
\frac{\rho(\xi(t, \bsw, \bsz))}{\varphi_0(\bsz)}
 -
\frac{\rho(\xi_h(t, \bsw, \bsz))}{\varphi_{0, h}(\bsz)}
\bigg| \bsrho(\bsw)\bsrho(\bsz) \, \rd \bsw\rd\bsz}_{I_{f, 1}} 
\nonumber\\
&+
\underbrace{\int_{\R^s}\int_{\R^s\setminus \Omega_f} \bigg|
\frac{\rho(\xi(t, \bsw, \bsz))}{\varphi_0(\bsz)}
 -
\frac{\rho(\xi_h(t, \bsw, \bsz))}{\varphi_{0, h}(\bsz)}
\bigg| \bsrho(\bsw)\bsrho(\bsz) \, \rd \bsw\rd\bsz}_{I_{f, 2}}
\,.
\end{align}

Rearranging the two fractions and using the triangle inequality, the first integral can be bounded by
\begin{align*}
I_{f, 1}
\leq & \underbrace{
\int_{\R^s}\int_{\Omega_f} \frac{|\rho(\xi(t, \bsw, \bsz)) - \rho(\xi_h(t, \bsw, \bsz))|}{\varphi_0(\bsz)}
\bsrho(\bsw) \bsrho(\bsz) \rd \bsw \rd \bsz
}_{I_{f, 3}}
\\
&+
\underbrace{\int_{\R^s}\int_{\Omega_f}
\frac{|\varphi_0(\bsz) - \varphi_{0, h}(\bsz)|\,\rho(\xi_h(t, \bsw, \bsz))}
{\varphi_0(\bsz) \varphi_{0, h}(\bsz)}
 \bsrho(\bsw) \bsrho(\bsz)\, \rd \bsw\rd\bsz}_{I_{f, 4}}.
\end{align*}

For the first term above, using the lower bound \eqref{eq:phi0-lower}
and the bound \eqref{eq:rho-diff} gives
\begin{align*}
I_{f, 3}
\leq \frac{1}{\sqrt{2\pi} \varphi_0(\bszero)}
\int_{\R^s}\int_{\Omega_f}
&K_0(\bsz)
\max \big\{ |\xi(t, \bsw, \bsz)|, |\xi_h(t, \bsw, \bsz)|\big\}
\\
&\cdot \big|\xi(t, \bsw, \bsz) - \xi_h(t, \bsw, \bsz)\big|
\bsrho(\bsw)\bsrho(\bsz)\,  \rd \bsw\rd\bsz.
\end{align*}
Expanding $\xi(t, \bsw, \bsz)$ using the equivalent of \eqref{eq:xi_h}, then using the triangle inequality along with
the lower bound \eqref{eq:phi0-lower} we can bound it from above by
\begin{align}
\label{eq:xi-bound}
|\xi(t, \bsw, \bsz)| \,&\leq  \frac{K_0(\bsz)}{ \varphi_0(\bszero)}
\bigg(|t| + |\phibar(\bsz)| + \sum_{i = 1}^s |w_i| |\varphi_i(\bsz)|\bigg)
\nonumber\\
&\leq \frac{K_0(\bsz)}{\varphi_0(\bszero)}
\bigg(|t| + \frac{\|\calG\|_{V^*} \big[ \|\ellbar\|_{V^*} + \sum_{i = 1}^s |w_i| \|\ell_i\|_{V^*}\big]}{\amin(\bsz)}\bigg),
\end{align}
where in the second inequality we have used the a priori bound \eqref{eq:LaxMil-G} for
$\phibar(\bsz) = \calG(\ubar(\cdot, \bsz))$ and $\varphi_i(\bsz) = \calG(u_i(\cdot, \bsz))$.
Since $V_h$ is conforming, the a priori bounds also hold for the FE approximations $\phibar_h(\bsz)$
and $\varphi_{i, h}(\bsz)$.
Then by using the lower bound on the FE approximation \eqref{eq:u_h-lower},
which holds for all $\bsz \in \Omega_h(\bsbeta)$ by Lemma~\ref{lem:Omega_h},
it follows that $\xi_h(t, \bsw, \bsz)$ also satisfies the bound \eqref{eq:xi-bound}, but with an extra factor 2.
Hence,
\[
\max \big\{ |\xi(t, \bsw, \bsz )|, |\xi_h(t, \bsw, \bsz)|\big\}
\leq \frac{2K_0(\bsz)}{\varphi_0(\bszero)}
\Bigg(|t| + \frac{\|\calG\|_{V^*} \big[ \|\ellbar\|_{V^*} + \sum_{i = 1}^s |w_i| \|\ell_i\|_{V^*}\big]}{\amin(\bsz)}\Bigg)\,.
\]

Substituting this along with \eqref{eq:xi-fe-err} into the bound on $I_{f, 3}$ gives
\begin{align}
\label{eq:I_fe3}
I_{f, 3}
&\leq\frac{2C_\calG\, h^{\tau' + \tau}}{\sqrt{2\pi} \varphi^3_0(\bszero)}
\int_{\R^s}\int_{\Omega_f} 
\frac{K_0^3(\bsz) \, \amax^4(\bsz) \|a(\cdot, \bsz)\|_{C^{0, 1}(\overline{D})}^4}
{\amin^9(\bsz)}
\Bigg( \|\ellbar\|_{L^2(D)} + \sum_{i = 1}^s |w_i| \|\ell_i\|_{L^2}
\nonumber\\&
+ \frac{2\|\ell_0\|_{L^2(D)} K_0(\bsz)}{\varphi_0(\bszero)}
\bigg[ |t| +  \frac{\|\calG\|_{V^*}\big[ \|\ellbar\|_{V^*} + \sum_{i = 1}^s |w_i| \|\ell_i\|_{V^*}\big]}{\amin(\bsz)}\bigg]
\Bigg)
\nonumber\\
&\cdot\bigg(|t| + \frac{\|\calG\|_{V^*} \big[ \|\ellbar\|_{V^*} + \sum_{i = 1}^s |w_i| \|\ell_i\|_{V^*}\big]}{\amin(\bsz)}\bigg)\bsrho(\bsw)\bsrho(\bsz) \, \rd \bsw\rd \bsz
\nonumber\\
&\leq\frac{2C_\calG\, h^{\tau' + \tau}}{\sqrt{2\pi} \varphi^3_0(\bszero)}
\int_{\R^s}\int_{\R^s}
\frac{K_0^3(\bsz)\, \amax^4(\bsz) \|a(\cdot, \bsz)\|_{C^{0, 1}(\overline{D})}^4}
{\amin^9(\bsz)}
\Bigg( \|\ellbar\|_{L^2(D)} + \sum_{i = 1}^s |w_i| \|\ell_i\|_{L^2}
\nonumber\\&
+ \frac{2\|\ell_0\|_{L^2(D)} K_0(\bsz)}{\varphi_0(\bszero)}
\bigg[ |t| +  \frac{\|\calG\|_{V^*}\big[ \|\ellbar\|_{V^*} + \sum_{i = 1}^s |w_i| \|\ell_i\|_{V^*}\big]}{\amin(\bsz)}\bigg]
\bigg)
\nonumber\\
&\cdot\bigg(|t| + \frac{\|\calG\|_{V^*} \big[ \|\ellbar\|_{V^*} + \sum_{i = 1}^s |w_i| \|\ell_i\|_{V^*}\big]}{\amin(\bsz)}\Bigg) \bsrho(\bsw)\bsrho(\bsz) \, \rd \bsw\rd \bsz
\nonumber\\
&=\, C_{f, 3}\,  h^{\tau' + \tau},
\end{align}
where the integral is finite by \eqref{eq:amax-bnd}--\eqref{eq:||a||-bnd}.

For the second term $I_{f, 4}$, we substitute in the FE error bound \eqref{eq:fe-err-L2} and
the lower bounds \eqref{eq:phi0-lower} and \eqref{eq:u_h-lower},
then use the property that the $\rho$ is bounded from above by $1/\sqrt{2\pi}$
to obtain
\begin{align}
\label{eq:I_fe4}
I_{f, 4}
&\leq h^{1 + \tau} \, \frac{2C_\calG\, \|\ell_0\|_{L^2(D)}}{\sqrt{2\pi} \varphi_0(\bszero)^2}
\int_{\R^s}\int_{\Omega_f} 
\frac{K_0^2(\bsz)\, \amax^4(\bsz) \|a(\cdot, \bsz)\|_{C^{0, 1}(\overline{D})}}
{\amin^9(\bsz)}
\bsrho(\bsw)\bsrho(\bsz) \,\rd \bsw \rd \bsz
\nonumber\\
&\leq h^{1 + \tau} \, \frac{2C_\calG\, \|\ell_0\|_{L^2(D)}}{\sqrt{2\pi} \varphi_0(\bszero)^2}
\int_{\R^s}\int_{\R^s}
\frac{K_0^2(\bsz)\, \amax^4(\bsz) \|a(\cdot, \bsz)\|_{C^{0, 1}(\overline{D})}}
{\amin^9(\bsz)}
\bsrho(\bsw)\bsrho(\bsz) \,\rd \bsw \rd \bsz
\nonumber\\
&= C_{f, 4}\  h^{1 + \tau},
\end{align}
where the integral is finite by \eqref{eq:amax-bnd}--\eqref{eq:||a||-bnd}.

Returning to \eqref{eq:f_h-err2}, the second integral  can be bounded using the
triangle inequality and the property that $\rho \leq 1/\sqrt{2\pi}$,
followed by the lower bounds \eqref{eq:phi0-lower} and \eqref{eq:phi_h-lower2}
(with $\varphi_0(\bszero) = \calG(u_0(\cdot, \bsz))$) to give
\begin{align*}
I_{f, 2} \leq &\frac{1}{\sqrt{2\pi}}\int_{\R^s} \int_{\R^s \setminus \Omega_f}
\bigg(\frac{1}{\varphi_0(\bsz)} + \frac{1}{\varphi_{0, h}(\bsz)}\bigg) \bsrho(\bsw)\bsrho(\bsz) \, \rd \bsw\rd \bsz
\\
\leq
&\frac{1}{\sqrt{2\pi}} \int_{\R^s} \bsrho(\bsw) \, \rd \bsw
\int_{\R^s \setminus \Omega_h(\bsbeta)} 
\bigg(\frac{K_0(\bsz)}{\varphi_0(\bszero)} + 
\frac{\amax(\bszero)}{\amin(\bszero)}\frac{2K_{0, h}(\bsz)}{\varphi_0(\bszero)}
\bigg) \bsrho(\bsz) \,\rd \bsz
\\
\leq &
\frac{1}{\sqrt{2\pi} \varphi_0( \bszero)}\bigg(1 + 2\frac{\amax(\bszero)}{\amin(\bszero)}\bigg)
\int_{\R^s \setminus \Omega_h(\bsbeta)} 
\big(K_0(\bsz) + K_{0, h}(\bsz)\big) \bsrho(\bsz) \,\rd \bsz.
\end{align*}
By the definitions \eqref{eq:K0} and \eqref{eq:K0h}, the term in the integral above can be simplified as 
\begin{align*}
&K_0(\bsz) +  K_{0, h}(\bsz)
\\
&= \amax(\bsz) \bigg( \frac{2}{\amin(\bszero)} + 2\frac{\|u_0(\cdot, \bszero)\|_{W^{1, \infty}(D)}}{\ell_{0, \inf}}
\sum_{j = 1}^\infty |z_j| \|a_j\|_{W^{1, \infty}(D)} + \frac{C_V h^{-d/2}}{\ell_{0, \inf}}\bigg)
\\
&\leq \exp\big(\|a_0\|_{L^\infty(D)}\big)\exp\Bigg(\sum_{j = 1}^s |z_j| \|a_j\|_{L^\infty(D)} \Bigg)
\exp\Bigg(\sum_{j = 1}^s |z_j| \|a_j\|_{W^{1, \infty}(D)}\Bigg)
\\
&\quad\cdot\bigg( \frac{2}{\amin(\bszero)} + \frac{2\|u_0(\cdot, \bszero)\|_{W^{1, \infty}(D)} + C_V h^{-d/2}}{\ell_{0, \inf}}\bigg)
\\
&= \exp\big(\|a_0\|_{L^\infty(D)}\big) \exp(\abs(\bsz) \cdot \bsbeta) 
\bigg( \frac{2}{\amin(\bszero)} 
+ \frac{2\|u_0(\cdot, \bszero)\|_{W^{1, \infty}(D)} + C_Vh^{-d/2}}{\ell_{0, \inf}}\bigg),
\end{align*}
where we have used the properties $\chi \leq \exp(\chi)$ and $\exp(\chi) \geq 1$ for $\chi > 0$,
the bound \eqref{eq:amax-bnd} and that, by definition, $\beta_j = \|a_j\|_{L^\infty(D)} + \|a_j\|_{W^{1, \infty}(D)}$.

Thus, $I_{f, 2}$ in \eqref{eq:f_h-err2} can be bounded by using \eqref{eq:int_not_Omega_h}
for $\Omega_h(\bsbeta)$ to give
\begin{align}
\label{eq:I_fe2}
I_{f, 2} \leq&
\frac{\exp\big(\|a_0\|_{L^\infty}\big)}{\sqrt{2\pi} \varphi_0(\bszero)}
\bigg(1 + 2 \frac{\amax(\bszero)}{\amin(\bszero)}\bigg)
\bigg( \frac{2}{\amin(\bszero)} + \frac{2\|u_0(\cdot, \bszero)\|_{W^{1, \infty}} + C_Vh^{-d/2}}{\ell_{0, \inf}}\bigg)
\nonumber\\
&\cdot\int_{\R^s \setminus \Omega_h(\bsbeta)} 
\exp(\abs(\bsz) \cdot \bsbeta) \bsrho(\bsz) \,\rd \bsz
\nonumber\\
\leq &
\frac{\exp\big(\|a_0\|_{L^\infty}\big)}{\sqrt{2\pi} \varphi_0(\bszero)}
\bigg(1 + 2\frac{\amax(\bszero)}{\amin(\bszero)}\bigg)
\bigg( \frac{2}{\amin(\bszero)} + \frac{2\|u_0(\cdot, \bszero)\|_{W^{1, \infty}} + C_V h^{-d/2}}{\ell_{0, \inf}}\bigg)
\nonumber\\&\cdot 
C_{\bsbeta, \tau} h^{d/2 + \tau + \tau'}
\nonumber\\
\leq & C_{f, 2} h^{\tau + \tau'}
\end{align}
where we have also used that $h$ is sufficiently small
to ensure that $C_{f, 2}$ is independent of $h$.

The final result \eqref{eq:f-fe-err} follows by using $I_{f, 1} \leq I_{f, 3} + I_{f, 4}$,
then substituting \eqref{eq:I_fe3}, \eqref{eq:I_fe4} and \eqref{eq:I_fe2} into \eqref{eq:f_h-err2}, 
noting that all constants are independent of $h$.
\end{proof}

\subsection{Combined QMC and FE error analysis}
Finally, we bound the total error of the combined FE and quadrature approximation
of the cdf and pdf using \eqref{eq:F_hN} and \eqref{eq:f_hN}, respectively,
for the specific case of a randomly shifted lattice rule \eqref{eq:QMC}
with $N$ a prime power.
Since the QMC rule is random, we now want to study the \emph{mean square error}
with respect to the random shift.

To bound the total mean square error of the approximation
$F_{h, N}(t)$, we first decompose it into the QMC and FE components,
as follows
\begin{equation}
\label{eq:err-decomp}
\bbE_\bsDelta\big[|F(t) - F_{h, N}(t)|^2\big] \leq
\bbE_\bsDelta\big[|F(t) - F_N(t)|^2\big]
+ \bbE_\bsDelta\big[|F_{N}(t) - F_{h, N}(t)|^2\big],
\end{equation}
and similarly for the error of $f_{h, N}(t)$. Here the subscript $\bsDelta$ signifies
that the expectation is with respect to the random shift $\bsDelta$ in the QMC rule.
Splitting the error in this way, the preintegration and QMC analysis is performed
only on the continuous problem, not the FE approximation.
Hence, there is no need to study the smoothness of the discrete approximation
with respect to the stochastic parameters.
For a lattice rule constructed by the component-by-component algorithm \cite{NK14}
the first term can be bounded by results from \cite{GKSr25} and then
the second term can be bounded similarly to the previous section.

\begin{theorem}
\label{thm:fe-err-cdfpdf}
Suppose that Assumption~\ref{asm:pde} holds. For all $\bsw \in \R^{s + 1}$,
let $\ell(\cdot, \bsw) \in H^{-1 + \tau'}(D)$ for some $\tau' \in (0, 1]$;
let $\calG \in H^{-1 + \tau}(D)$ for some $\tau \in [0, 1]$ be a positive linear functional;
and 
for $\eps \in (0, 1/2]$, let $\psi$ be a Gaussian weight function with $\mu \in (0, \eps)$.
For $h > 0$ sufficiently small, $N \in \N$ a prime power and $t \in [t_0, t_1]$, 
the combined FE and QMC approximations $F_{h, N}(t)$ \eqref{eq:F_hN}
and $f_{h, N}(t)$ \eqref{eq:f_hN} of the cdf
and pdf of the quantity of interest $X = \calG(u(\cdot, \bsy))$ satisfy
\begin{align}
\label{eq:F-err}
\sqrt{\bbE_\bsDelta\big[|F(t) - F_{h, N}(t)|^2\big]}
&\leq C_{F, t, \tau, \eps} \, \big( h^{\tau' + \tau} + N^{-1 + \eps}\big),
\quad \text{and}\\
\label{eq:f-err}
\sqrt{\bbE_\bsDelta\big[|f(t) - f_{h, N}(t)|^2\big]}
&\leq C_{f, t, \tau, \eps} \, \big(h^{\tau' + \tau} + N^{-1 + \eps}\big),
\end{align}
where $C_{F, t, \tau, \eps}, C_{f, t, \tau, \eps} < \infty$ are independent of $h$ and $N$, but depend on $s$.
\end{theorem}

\begin{proof}
The first term in the decomposition \eqref{eq:err-decomp} can be bounded by \cite[Theorem 6.1]{GKSr25},
\begin{equation}
\label{eq:qmc-err}
\bbE_\bsDelta\big[|F(t) - F_{h, N}(t)|^2\big] \leq \big[ C_{\bsgamma, s, \eps} N^{-1 + \eps} \big]^2,
\end{equation}
with $C_{\bsgamma, s, \eps}$ independent of $N$, but depending on $s$ and the weights $\bsgamma$
as specified in the proof of \cite[Theorem 6.1]{GKSr25}.

For the second term in \eqref{eq:err-decomp}, using the simplification \eqref{eq:F_hN}, along with
the linearity of a QMC rule \eqref{eq:QMC}, we can write the mean square FE error as
\begin{align*}
\bbE_\bsDelta\big[|F_N(t) - F_{h, N}(t)|^2\big]
\,&=\, \bbE_\bsDelta\big[|Q_{2s, N}\big((\Phi(\xi(t, \cdot)) - \Phi(\xi_h(t, \cdot))\big)|^2\big]
\\
&=\,
\bbE_\bsDelta\Bigg[\frac{1}{N^2}\bigg|\sum_{k = 0}^{N - 1}
\big[\Phi(\xi(t, \bsq_k)) - \Phi(\xi_h(t, \bsq_k))\big]
\bigg|^2\Bigg].
\end{align*}
By the Cauchy--Schwarz inequality for sums we can bound this by
\begin{align*}
&\bbE_\bsDelta\big[|F_N(t) - F_{h, N}(t)|^2\big]
\\
&\leq
\bbE_\bsDelta\Bigg[\frac{1}{N^2} \Bigg(\sum_{k = 0}^{N - 1}
|\Phi(\xi(t, \bsq_k)) - \Phi(\xi_h(t, \bs1_k))|^2\Bigg)
\Bigg( \sum_{n = 0}^{N - 1} 1^2\Bigg)
\Bigg]
\\
&=\, \bbE_\bsDelta\Bigg[\frac{1}{N} \sum_{n = 0}^{N - 1}
|\Phi(\xi(t, \tauhat_n(\bsDelta))) - \Phi(\xi_h(t, \tauhat_n(\bsDelta)))|^2\Bigg]
\\
&=\, \bbE\big[|\Phi(\xi(t, \cdot)) - \Phi(\xi_h(t, \cdot))|^2\big],
\end{align*}
where we have used the property that randomly shifted QMC rules are unbiased.
Thus, we can write the mean square FE error for the cdf as
\begin{equation*}
\bbE_\bsDelta\big[|F_N(t) - F_{h, N}(t)|^2\big]
\leq
\int_{\R^{2s}} \big|\Phi(\xi(t, \bsy))
- \Phi(\xi_h(t, \bsy))\big|^2
\bsrho(\bsy) \, \rd \bsy.
\end{equation*}
The right hand side above is similar to \eqref{eq:F_h-err}, but with a power of 2 on the integrand,
and thus, it can be bounded in the same way. To account for the power of 2, the only change is to also truncate the domain of integration for $\bsw$ to $\Omega_h(\bszero)$.\footnote{Alternatively, one could stick with truncating only the $\bsz$ domain, but use a larger box by replacing the power of $h$ in \eqref{eq:kappa} by $d + 2(\tau + \tau')$.}
Then, following the same steps as in the proof of Theorem~\ref{thm:fe-err} 
we can show the bound
\begin{equation}
\label{eq:F_h-err1}
\bbE_\bsDelta\big[|F_N(t) - F_{h, N}(t)|^2\big]
\leq \widetilde{C}_{F, t, \tau} h^{2(\tau + \tau')},
\end{equation}
where $\widetilde{C}_{F, t,\tau}$ is again independent of $h$ and $\bsy$. 

Substituting \eqref{eq:qmc-err} and \eqref{eq:F_h-err1} into \eqref{eq:err-decomp} then taking the square root
gives the required result \eqref{eq:F-err}.
The proof of \eqref{eq:f-err} follows in the same way.
\end{proof}

To emphasise the significance of this result, 
note that the convergence above for the combined error in computing the cdf and pdf
is the same rate, in both $h$ and $N$, as for the much simpler problem of computing the expected 
value of the QoI, as in, e.g., \cite{GrKNSSS15}.
Note that the constants are not the same, with the constant here possibly depending on $s$. 

\section{Conclusion}
This paper analysed the FE component of density estimation
using preintegration for random elliptic PDEs of the form \eqref{eq:pde}.
In particular, it was proven that the FE approximation of the QoI
also satisfies the important monotonicity condition
and the analysis of the FE error for the cdf and pdf presented here
provides a complete analysis of the error for the method from the earlier paper \cite{GKSr25}.
Importantly, the convergence rates shown are the same as for 
computing the expected value.

\bibliographystyle{plain}

\end{document}